\documentclass[12pt,reqno]{amsart}

\usepackage[T1]{fontenc} 
\usepackage{amssymb}
\usepackage{amscd}
\usepackage{upref}

\theoremstyle{plain}
\newtheorem{theorem}{Theorem}[section]
\newtheorem{lemma}[theorem]{Lemma}
\newtheorem{proposition}[theorem]{Proposition}

\theoremstyle{definition}

\theoremstyle{remark}
\newtheorem{remark}[theorem]{Remark}
\newtheorem*{acknowledgements}{Acknowledgements}

\numberwithin{equation}{section}

\newcommand{\abs}[1]{\lvert #1\rvert}
\DeclareMathOperator{\Ad}{Ad}
\DeclareMathOperator{\ad}{ad}

\DeclareMathOperator{\cohom}{cohom}
\DeclareMathOperator{\CP}{\mathbb CP}

\DeclareMathOperator{\End}{End}

\DeclareMathOperator{\diag}{diag}
\DeclareMathOperator{\GL}{\textsl{GL}}
\DeclareMathOperator{\gl}{\mathfrak{gl}}
\DeclareMathOperator{\Gr}{Gr}
\DeclareMathOperator{\grad}{grad}
\DeclareMathOperator{\Gro}{\widetilde{Gr}}

\DeclareMathOperator{\HP}{\mathbb HP}

\DeclareMathOperator{\im}{Im}
\DeclareMathOperator{\Isom}{Isom}
\newcommand{\inp}[3][]{\left\langle #2,#3\right\rangle_{#1}}

\DeclareMathOperator{\LieE}{\textsl E}
\DeclareMathOperator{\lieE}{\mathfrak e}
\DeclareMathOperator{\LieF}{\textsl F}
\DeclareMathOperator{\lieF}{\mathfrak f}
\DeclareMathOperator{\LieG}{\textsl G}
\DeclareMathOperator{\lieG}{\mathfrak g}
\DeclareMathOperator{\LSP}{\mathfrak{sp}}

\newcommand{\norm}[2][]{\left\lVert #2\right\rVert_{#1}}

\DeclareMathOperator{\RP}{\mathbb RP}

\DeclareMathOperator{\SL}{\textsl{SL}}
\DeclareMathOperator{\Sl}{\mathfrak{sl}}
\DeclareMathOperator{\SO}{\textsl{SO}}
\DeclareMathOperator{\so}{\mathfrak{so}}
\DeclareMathOperator{\SP}{\textsl{Sp}}
\DeclareMathOperator{\Special}{\textsl S}

\DeclareMathOperator{\Spin}{\textsl{Spin}}
\DeclareMathOperator{\spin}{\mathfrak{spin}}
\DeclareMathOperator{\stab}{stab}

\DeclareMathOperator{\SU}{\textsl{SU}}
\DeclareMathOperator{\su}{\mathfrak{su}}
\DeclareMathOperator{\Sym}{Sym}

\DeclareMathOperator{\UM}{\mathcal U}
\DeclareMathOperator{\Un}{\textsl U}
\DeclareMathOperator{\un}{\mathfrak u}

\makeatletter

\newcommand{\Dyo}{\mathord{\circ}}
\newcommand{\Dyx}{\mathord{\times}}

\newcommand{\ltwoa}{\mathord{\mathpalette\tw@a<}}
\newcommand{\rtwoa}{\mathord{\mathpalette\tw@a>}}
\newcommand{\tw@a}[2]{\ooalign{\hfil$#1 #2$\hfil\crcr
$#1 \Relbar\joinrel\Relbar$\crcr}}

\newcommand{\lthreea}{\mathord{\mathpalette\thr@a<}}
\newcommand{\rthreea}{\mathord{\mathpalette\thr@a>}}
\newcommand{\thr@a}[2]{\ooalign{\hfil$#1 #2$\hfil\crcr
\raise.3ex\hbox{$#1 \relbar\mathrel{\mkern-4mu}\relbar$}\crcr
$#1 \relbar\mathrel{\mkern-4mu}\relbar$\crcr
\lower.3ex\hbox{$#1 \relbar\mathrel{\mkern-4mu}\relbar$}\crcr}}

\newcommand{\Dyatop}[2]{\mathord{\mathpalette{\Dy@top{#1}{#2}}{}}}
\newcommand{\Dy@top}[3]{\overset{#3 #1}{#3 #2}}

\makeatother

\def\englishdate#1/#2/#3.{#2/#1/#3.}   

\begin{document}
\title{Quaternionic K\"ahler Manifolds of Cohomogeneity One}

\author{Andrew Dancer}
\address[Dancer]{Department of Mathematics and Statistics\\
McMaster University\\
Hamilton\\
Ontario L8S 4K1\\
Canada}
\email{dancer@icarus.math.mcmaster.ca}

\author{Andrew Swann}
\address[Swann]{Department of Mathematical Sciences\\
University of Bath\\
Claverton Down\\ 
Bath BA2 7AY\\
England}
\email{A.F.Swann@maths.bath.ac.uk}


\begin{abstract}
  Classification results are given for (i)~compact quaternionic K\"ahler
  manifolds with a cohomogeneity-one action of a semi-simple group,
  (ii)~certain complete hyperK\"ahler manifolds with a
  cohomogeneity\kern0pt-two action of a semi-simple group preserving each
  complex structure, (iii)~compact $3$-Sasakian manifolds which are
  cohomogeneity one with respect to a group of $3$-Sasakian symmetries.
  Information is also obtained about non-compact quaternionic K\"ahler
  manifolds of cohomogeneity one and the cohomogeneity of adjoint orbits in
  complex semi-simple Lie algebras.
\end{abstract}

\maketitle

\tableofcontents

\newpage
\section{Introduction}
\label{sec:introduction}
A Riemannian manifold $(M,g)$ is said to be of cohomogeneity one with
respect to a group~$G$ if $G$~acts isometrically with generic orbit of real
codimension one.  There has been considerable interest in studying Einstein
manifolds of this type, because the Einstein condition then reduces to a
set of ordinary differential equations.

Our aim is to investigate quaternionic K\"ahler manifolds which are of
cohomogeneity one.  Recall that a manifold of dimension~$4n$ is
quaternionic K\"ahler if the holonomy reduces to a subgroup of $\SP(n)\SP(1)$
(see \cite{Besse:Einstein}, for example).  Such manifolds are always
Einstein.  Each quaternionic K\"ahler manifold~$M$ has associated to it a
twistor space~$Z$: a complex contact manifold of complex dimension~$2n+1$
mapping onto~$M$ with fibers rational curves of normal bundle~$2n\mathcal
O(1)$.

Quaternionic isometries of~$M$ give rise to holomorphic symmetries of~$Z$
preserving the complex contact structure.  If $M$ has cohomogeneity one,
then, in many cases, moment map techniques show that $Z$ is related to a
coadjoint orbit in a complex Lie algebra.  We are thus lead to discuss the
action of a real group~$G$ on $G^{\mathbb C}$-orbits in the complex Lie
algebra $\mathfrak g^{\mathbb C}$.  When $G$~is compact, we prove a series
of monotonicity results that say that the cohomogeneity of such orbits
varies in a natural way with respect to various partial orders on the
orbits.  This enables us to narrow down the possibilities for $M$ and~$Z$,
and in the compact case, we are able to determine all examples~$(M,g,G)$
when $G$~is semi-simple.

In the compact case, partial results were obtained by Alekseevsky \&\ 
Podest\`a~\cite{Alekseevsky-Podesta:co1}.  Most of our techniques are
different, though we do use one of their arguments.  In our paper, the
geometry of nilpotent orbits plays a notable r\^ole and results of
Brylinski \&\ Kostant \cite{Brylinski-Kostant:nilpotent} on shared orbits
and of Beauville~\cite{Beauville:Fano} on Fano contact manifolds are
applied.

It turns out that the compact quaternionic K\"ahler manifolds that we
obtain are all symmetric spaces.  Kollross~\cite{Kollross:thesis} has
recently classified cohomogeneity\kern0pt-one actions of reductive groups
on compact symmetric spaces and our results are consistent with his.
However, one difference to note is that he considers two such actions to be
equivalent if they have the same orbits.  Thus in some cases he gives a
space with an action of a reductive group, where we obtain the same space
with a semi-simple symmetry group.  One of the main open questions in
quaternionic K\"ahler geometry is whether there exist any compact
non-symmetric examples with positive scalar curvature.  Our work shows that
there are no such examples of cohomogeneity one with respect to a
semi-simple group.

Our techniques, using moment maps and cohomogeneities of adjoint orbits,
also enable us to tackle classification problems for two other types of
geometric structure: hyperK\"ahler and $3$-Sasakian.
In~\cite{Bielawski:homogeneous}, Bielawski showed how a certain class of
complete hyperK\"ahler manifolds always arises as coadjoint orbits.  We
analyse this situation in more detail to classify those of cohomogeneity
two.

There has been much recent interest in $3$-Sasakian manifolds because they
provide new examples of compact Einstein manifolds with positive scalar
curvature.  As shown in~\cite{Boyer-GM:three-Sasakian}, there is a close
relationship between $3$-Sasakian manifolds, hyperK\"ahler structures and
quaternionic K\"ahler orbifolds.  Using knowledge of cohomogeneities of
adjoint orbits and extending the results of
Beauville~\cite{Beauville:Fano}, we are able to classify $3$-Sasakian
manifolds whose group of $3$-Sasakian symmetries acts with cohomogeneity
one.

\begin{acknowledgements}
  The first author is partially supported by \textsc{nserc} grant
  \textsc{opg}\oldstylenums{0184235}.  He thanks Claude LeBrun for useful
  discussions.  The second author would like to thank the organisers of the
  conference on Complex and Symplectic Geometry in Cortona for a
  stimulating atmosphere and is grateful to Francis Burstall, Alastair
  King, Piotr Kobak and Kris Galicki for useful discussions.  The second
  author is partially supported by the \textsc{epsrc} of Great Britain.  We
  thank Roger Bielawski for comments on an earlier version of this
  manuscript.
\end{acknowledgements}

\section{Hypotheses}
\label{sec:hypotheses}
Let $M$~be a quaternionic K\"ahler manifold of dimension~$4n$, with twistor
space~$\pi\colon Z\to M$, see~\cite{Besse:Einstein}.  Suppose that $G$~is a
compact semi-simple Lie group that acts quaternionically on~$M$ with
cohomogeneity one.  Without loss of generality we may also assume that
$G$~is connected.  Write~$\mathfrak g$ for the Lie algebra of~$G$ and
$\mathfrak g^{\mathbb C}$ for the complexification of~$\mathfrak g$.  Let
$M_0$ be the union of principal orbits in~$M$.  General information about
$G$-manifolds may be found in~\cite{Bredon:transformation}.

\section{The Action on the Twistor Space}
\label{sec:twistor}
The $G$-action lifts to a group of holomorphic contact transformations
on~$Z$.  We then have a map of sheaves
\begin{equation*}
  \alpha\colon \mathcal O\otimes\mathfrak g^{\mathbb C} \longrightarrow TZ
\end{equation*}
given by mapping elements of~$\mathfrak g^{\mathbb C}$ to the corresponding
vector fields on~$Z$ (cf.~\cite{Hitchin:isomonodromic}).

\begin{lemma}
  At each point~$z$ of~$Z$, the image of~$\alpha$ has real dimension either
  $4n$ or $4n+2$.
\end{lemma}

\begin{proof}
  As the action of~$G$ is of cohomogeneity one on~$M$, the lifted action
  has cohomogeneity at most three.  Thus the real dimension of the image
  of~$\alpha$ is at least $4n-1$.  However, the image of~$\alpha$ is
  complex, so the real dimension is even, and thus either $4n$ or $4n+2$.
\end{proof}

\begin{lemma}
  \label{lem:U0}
  Let $U_0$ be the set of points~$z\in Z$ where the dimension of the image
  of~$\alpha$ is~$4n+2$.  Then $U_0$ is open and either empty or dense.
\end{lemma}

\begin{proof}
  Taking exterior powers, we get a bundle map
  \begin{equation*}
    \Lambda^{2n+1}\alpha \colon \mathcal O \otimes \Lambda^{2n+1} \mathfrak
    g^{\mathbb C} \to K_Z^{-1},
  \end{equation*}
  where $K_Z=\Lambda^{2n+1}T^*Z$ is the canonical bundle.  The set~$U_0$ is
  just the complement of the zero set of the analytic
  map~$\Lambda^{2n+1}\alpha$.  Thus either $U_0$~is empty or $U_0$~is open
  and dense.
\end{proof}

We thus have two possibilities: either $U_0$~is empty or it is not.
Consider the complexified group~$G^{\mathbb C}$.  This need not act on~$Z$
for reasons of completeness, but it does act in the sense of groupoids.
The case $U_0\ne\varnothing$ corresponds exactly to the existence of an
open orbit for this action of~$G^{\mathbb C}$ on~$Z$.  The treatment of
this case will start in~\S\ref{sec:moment}.

\section{Twistor Spaces with No Open Orbits}
\label{sec:no-open}
Let us consider the case when the twistor space has no open $G^{\mathbb
C}$-orbits.  This is, in fact, the simplest case.  Let $\widehat M$ be
the universal cover of the union of principal orbits in~$M$.  We first show
that $\widehat M$ admits a hypercomplex structure.  This structure can
taken to be $G$-invariant and the aim is then to show that $\widehat
M$~fibres over a homogeneous quaternionic manifold.  The assumption that
$G$~is compact implies that the quotient is a Wolf space and the structure
on~$M$ can then be analysed metrically.  We are then able to show that the
only compact example is given by the action of~$\SP(n)$ on~$\HP(n)$.

\begin{proposition}
  If $U_0$ is empty, then $M$~is locally hypercomplex and each point
  of~$\widehat M$ admits a $\widetilde G$-invariant neighbourhood with a
  $\widetilde G$-invariant hypercomplex structure, where $\widetilde G$ is
  the universal cover of~$G$.
\end{proposition}

\begin{proof}
  The complex rank of~$\alpha$ is~$2n$.  The definition of~$\alpha$ implies
  that its image is integrable in the sense of Frobenius, so we have a
  foliation of $Z_0$ by leaves of complex dimension $2n$ tangent to the
  image of $\alpha$.  If we compose $\alpha$~with projection to the
  horizontal distribution on~$Z_0$, then we get a map~$\alpha^{\mathcal
  H}$, whose real rank is at least $4n-1$ as the action on~$M$ has
  cohomogeneity one.  However, the image of~$\alpha^{\mathcal H}$ is
  complex, so $\alpha^{\mathcal H}$ is a surjection, and the leaves of the
  foliation are transverse to twistor lines.  This means that locally we
  have a holomorphic projection $Z_0\to\CP(1)$ and $M$~is locally
  hypercomplex.

  The map~$\alpha$ is $G$-equivariant, so the foliation determined by the
  image of~$\alpha$ is invariant under the infinitesimal action of~$G$.
  However, we have assumed that $G$~is connected, so $G$~preserves the
  foliation and the local hypercomplex structure.
  
  Now consider $\widehat M$.  Topologically, $\widehat M$~is a product
  $\mathbf I\times \widetilde G/H$ for some interval~$\mathbf I$.  As
  $\widetilde G$~is connected and $\widehat M$~is simply-connected,
  the exact homotopy sequence for a fibration implies that $H$ is also
  connected.  Fix $t_0$ in~$\mathbf I$.  There is a neighbourhood~$V$
  of $x=t_0\times H$ in~$\widehat M$ on which $\widehat M$ has a
  hypercomplex structure $I,J,K$.  Consider the $\widetilde G$-orbit
  of~$I$.  This meets the twistor line~$\widehat Z_x$ in a set which
  has empty interior and is an $H$-orbit.  As $H$ is connected, this
  orbit is a single point.  Thus we have a $\widetilde G$-invariant
  hypercomplex structure on $\mathbf J\times \widetilde G/H$, where
  $\mathbf J\subset\mathbf I$ is an open interval containing~$t_0$.
\end{proof}

We say a quaternionic manifold has a compatible hypercomplex structure, if
there exist global integrable sections $I$, $J$ and $K$ of~$\mathcal G$
satisfying the quaternion identities.

\begin{proposition}
  \label{prop:sp1}
  Suppose $M$~is a quaternionic K\"ahler manifold of positive scalar
  curvature.  If $M$~has a compatible hypercomplex structure~$I,J,K$,
  then $M$~admits an infinitesimal action of~$\SP(1)\times\mathbb R$ with
  the following properties:
  \begin{enumerate}
  \item[(i)] $\SP(1)$~acts isometrically;
  \item[(ii)] If $V$~is the vector field generating the action of~$\mathbb
    R$, then the action of~$\SP(1)$ is generated by~$IV$, $JV$ and $KV$;
  \item[(iii)] $V$~preserves $I$, $J$ and~$K$;
  \item[(iv)] $L_{IV}I=0$ and $L_{IV}J=K$, and these formul\ae{} remain
    valid for any cyclic permutation of~$I,J,K$.
  \end{enumerate}
\end{proposition}

\begin{proof}
  The construction follows~\cite[\S5]{Swann:MathAnn}.  Let $\nabla$~be the
  Levi-Civita connection and let $\nabla'$~be the Obata connection.  The
  twistor operator~$D$ is the composition
  \begin{equation*}
    \begin{CD}
      H
      @>\nabla''>>
      H\otimes EH \cong ES^2H \oplus E
      @>p>>
      ES^2H,
    \end{CD}
  \end{equation*}
  where $\nabla''$~is any quaternionic connection and 
  $p$~is projection.  The operator~$D$ is
  independent of~$\nabla''$.  The hypercomplex structure defines a
  section~$h$ of~$H$ such that $\nabla'h=0$.  In particular, $Dh=0$.  Thus
  $\nabla h$~lies in the module~$E$.  Let $e=\nabla h$.  By
  \cite[Lemma~5.7]{Swann:MathAnn}, $\nabla e=\lambda h$, for some
  constant~$\lambda$ (a positive constant times the scalar
  curvature of~$M$).
  
  Let $V=e\tilde h-\tilde e h$, where $\tilde e=je$ and $\tilde h=jh$.
  Then $IV=eh+\tilde e\tilde h$, $JV=i(e\tilde h+\tilde eh)$ and
  $KV=i(eh-\tilde e\tilde h)$.  We claim that $V$~is quaternionic and $IV$,
  $JV$ and $KV$ are Killing.  Note that
  \begin{equation*}
    \begin{split}
      T^*M\otimes T^*M
      &\cong S^2T^*M + \Lambda^2T^*M\\
      &\cong (\mathbb R+\Lambda^2_0E+S^2ES^2H)\\
      &\qquad + (S^2E+S^2H+\Lambda^2_0ES^2H).
    \end{split}
  \end{equation*}
  Killing vector fields are characterised by having covariant derivative
  in~$\Lambda^2T^*M$ and quaternionic vector fields are those with
  covariant derivative in
  \begin{equation*}
    \gl(n,\mathbb H)+\LSP(1)\cong \mathbb R + \Lambda^2_0E + S^2E + S^2H.
  \end{equation*}
  Now 
  \begin{equation}
    \nabla V=\lambda h\wedge\tilde h+e\wedge\tilde e
    \label{eq:nabla-V}
  \end{equation}
  which lies in $\mathbb R+\Lambda^2_0E\subset \gl(n,\mathbb H)$ so $V$~is
  not only quaternionic but it also preserves the hypercomplex structure.
  Similarly computation of the covariant derivatives of $IV$, $JV$ and $KV$
  shows that these three vector fields are Killing,
  \begin{equation*}
    [IV,JV]=2(\norm e^2+\lambda\norm h^2)KV
  \end{equation*}
  and
  \begin{equation*}
    L_{IV}J = 2(\norm e^2+\lambda\norm h^2) K.
  \end{equation*}
  Here $\norm e^2=e\wedge\tilde e$ and $\norm h^2=h\wedge \tilde h$, which
  are identified with functions via the symplectic forms on $E$ and $H$.
  Now $\nabla(h\wedge \tilde h)=0$ and $\nabla(e\wedge\tilde e)=0$, so
  $\norm e^2+\lambda\norm h^2$ is a constant.  Thus dividing~$V$ by the
  constant $2(\norm e^2+\lambda\norm h^2)$ yields a vector field with the
  required properties.
\end{proof}

\begin{remark}
  If the scalar curvature of~$M$ is negative then the above proof goes
  through unchanged provided~$\lambda\norm h^2\ne -\norm e^2$.  If
  $\lambda\norm h^2=-\norm e^2$ then we still get $V$ but now $IV$, $JV$
  and $KV$ commute.  We suspect this case does not arise.
\end{remark}

If $M$~is of cohomogeneity one with an invariant hypercomplex structure,
then we may choose the section~$h$ to be $G$-invariant.  This leads to:

\begin{proposition}
  \label{prop:associated}
  If $U_0$~is empty and $M$~has positive scalar curvature, then
  $\widehat M$~is the associated bundle of a compact Wolf space.
\end{proposition}

\begin{proof}
  From the previous proof we have that $\nabla V$~lies in~$S^2T^*\widehat
  M$.  As $\widehat M$ is simply connected, this implies that the vector
  field~$V$ is a gradient, say $V=\grad \rho$.  Now $h$~is $\widetilde
  G$-invariant, so its covariant derivative~$e$ is also and hence the
  $\widetilde G$-action commutes with the (infinitesimal) action
  of~$\SP(1)\times\mathbb R$ on~$\widehat M$.
  
  We claim that the function~$\rho$ is $\SP(1)$-invariant.  As $IV$
  commutes with $V$, we have $d((IV)\rho)=L_{IV}V^\flat=0$ and hence
  $(IV)\rho$~is constant.  Similarly, $(JV)\rho$ and $(KV)\rho$ are
  constant.  This implies $[JV,KV]\rho=(JV)((KV)\rho)-(KV)((JV)\rho)=0$.
  However, $[JV,KV]$~is a non-zero constant times~$IV$, hence $\rho$~is
  $\SP(1)$-invariant.  Similarly, $\rho$~is $\widetilde G$-invariant, as
  $\widetilde G$~is semi-simple and so each element of its Lie algebra is a
  sum of commutators.
  
  As the action of $\widetilde G$ on $\widehat M$ is cohomogeneity one,
  each component of a generic level set of $\rho$ will be an orbit of
  $\widetilde G$ and hence compact. The infinitesimal action of $\SP(1)$
  therefore integrates to a genuine group action, and the quotient of a
  component of a generic level set by $\SP(1)$ will be $\widetilde
  G$-homogeneous and hence smooth.
  
  The proof of \cite[Theorem~5.1]{Swann:MathAnn} may now be modified to
  show that such a quotient is a quaternionic K\"ahler manifold~$N$.  (The
  two changes required are to replace $\mu$ by~$\rho$ and to note that the
  covariant derivative of~$\omega_I$ is now a linear combination
  of~$\omega_J$ and $\omega_K$.)
  
  By~\cite{Alekseevsky:compact,Alekseevsky:transitive}, $N$~is a Wolf
  space.  By~\cite{Pedersen-PS:hypercomplex} the topology of~$\widehat M$
  is that of~$\UM(N)$.
\end{proof}

\begin{theorem}
  If $U_0$~is empty \textup(i.e., $Z$ has no open $G^{\mathbb
  C}$-orbit\textup) and $M$~has positive scalar curvature, then the metric
  on~$\widehat M$ is a member of the one-parameter family of quaternionic
  K\"ahler metrics given in~\cite[Theorem~3.5]{Swann:MathAnn} on the
  associated bundle of a compact Wolf space.
\end{theorem}

\begin{proof}
  The principal orbit of~$\UM(N)$ is~$G/H$, where $N=G/H\SP(1)$ is the Wolf
  space.  The Lie algebra~$\mathfrak g$ of~$G$ splits under the action
  of~$\Ad(H)$ as
  \begin{equation*}
    \mathfrak g = \mathfrak h \oplus \mathfrak p_0 \oplus \mathfrak p_1,
  \end{equation*}
  where $\mathfrak p_0$~is a trivial module and $\mathfrak p_1$~is a sum of
  non-trivial $\Ad(H)$-modules.  As $N$~is a Wolf space, $\mathfrak p_0$~is
  three-dimensional and isomorphic to~$\LSP(1)$, and $\mathfrak p_1$~is
  isomorphic to the tangent space of~$N$.
  
  We may identify~$\widehat M=\UM(N)$ with $(t_0, t_1) \times G/H$ for some
  subinterval $(t_0 , t_1)$ of $\mathbb R$, in such a way that the metric
  is written as $g=dt^2+g_t$, where $g_t$~is a homogeneous metric on~$G/H$
  for each $t$.  By Schur's Lemma, $\mathfrak p_0$ and $\mathfrak p_1$ are
  orthogonal for~$g_t$ and the hyperHermitian structure preserves the
  splitting
  \begin{equation*}
    T\UM(N)=\left(\left\langle\frac\partial{\partial t}\right\rangle \oplus
      \mathfrak p_0\right) \oplus \mathfrak p_1,
  \end{equation*}
  Thus, $\mathfrak p_1$~is the orthogonal complement to the quaternionic
  span of~$V$.  Equation~\eqref{eq:nabla-V} shows that on~$\mathfrak p_1$,
  the vector field~$V$ acts conformally.  Moreover, $\SP(1)$~acts
  isometrically and $\mathfrak p_1$~is irreducible as a representation
  of~$H\times\SP(1)$.  Thus the restriction of~$g_t$ to~$\mathfrak p_1$ is
  a multiple of the metric on~$TN$, which is itself a constant
  times~$\inp[\mathfrak p_1]\cdot\cdot$, the restriction to~$\mathfrak p_1$
  of minus the Killing form on~$\mathfrak g$.  Let $f(t)$ be the function
  such that $g_t|_{\mathfrak p_1} = f(t)^2 \inp[\mathfrak p_1]\cdot\cdot$.
  
  The restriction of~$g$ to~$\left\langle \frac\partial{\partial t}
  \right\rangle + \mathfrak p_0$ defines a Bianchi~IX hyperHermitian
  metric, so as in~\cite{Hitchin:isomonodromic} we can diagonalise
  this part of the metric for all~$t$.  More concretely, we can find a
  basis $\{X_1,X_2,X_3\}$ for $\mathfrak p_0 \cong \LSP(1)$ satisfying
  $[X_1,X_2]=2X_3$ etc., such that $g(X_a,X_b) = \delta_{ab} h_a^2$
  for $a,b=1,2,3$, where $h_a$ are functions of~$t$.  Note that by
  Proposition~\ref{prop:associated}, the hyperHermitian structure
  on~$\widehat M$ induces the quaternionic structure on the Wolf space, and
  in particular, the complex structures act via $\ad X_a$ on~$\mathfrak
  p_1$.
  
  Writing $J_1,J_2,J_3$ for $I,J,K$, we may argue
  as in~\cite{Dancer-Swann:hK-cohom1} and derive the formula
  \begin{equation*}
    J_a|_{\mathfrak p_1} = \frac{h_a}{2 f f'} \ad X_a \qquad\text{for
    $a=1,2,3$.}
  \end{equation*}
  The only change needed in the proof, is the contribution to
  $\nabla_{Y_i}(J_a\partial/\partial t)$ from terms $(\nabla
  J_a)\partial/\partial t$.  But as $\widehat M$~is quaternionic K\"ahler,
  each of these terms lies in~$\mathfrak p_0$ and so does not contribute
  to~$J_a$ on~$\mathfrak p_1$.
  
  As $J_a^2=-1$, we see that the functions $h_a^2$ are all equal.  This
  means that the metric on~$\widehat M=\UM(N)$ is of the form considered
  in~\cite[\S3]{Swann:MathAnn}.  Lemma~3.4 in~\cite{Swann:MathAnn}, implies
  that $g$~is of the form in~\cite[Theorem~3.5]{Swann:MathAnn}.
\end{proof}

\begin{theorem}
  \label{thm:hc-compact}
  Suppose $M$ is a compact quaternionic K\"ahler manifold of positive scalar
  curvature with a cohomogeneity-one action of a compact semi-simple
  group~$G$.  If $G^{\mathbb C}$~has no open orbit on the twistor space
  of~$M$, then $M=\HP(n)$ with its symmetric metric and $G$ is the
  subgroup~$\SP(n)$ of the full isometry group~$\SP(n+1)$.
\end{theorem}

\begin{proof}
  Without loss of generality we may assume that $G$~is connected.  As
  $M$~is simply-connected~\cite{Salamon:Invent}, the exact homotopy
  sequence implies that the compact one-dimensional manifold~$M/G$ is a
  closed interval.  Thus we have two special orbits $G/H_1$ and $G/H_2$.
  If the principal orbit is~$G/H$, then the fact that $M$~is a smooth
  manifold implies that $H_i/H$~is a sphere for $i=1,2$.
  
  The form of the metric in~\cite{Swann:MathAnn} implies that one of the
  special orbits is a point and the other is the underlying Wolf space~$W$.
  Thus $G/H$~is a sphere~$S^{4n-1}$ which implies $G=\SP(n)$ and
  $W=\HP(n-1)$, giving $M=\HP(n)$.  However, $\HP(n)$~only admits one
  hyperK\"ahler metric~\cite{Salamon:Invent}.
\end{proof}

\section{Moment Maps}
\label{sec:moment}
We assume from now on that the groupoid action of~$G^{\mathbb C}$ has an
open orbit~$U_1$ on the twistor space~$Z$, or equivalently that the
set~$U_0$ of Lemma~\ref{lem:U0} is non-empty.  The idea is to use a moment
map construction to identify an open set of~$Z$ with an open set in a
projectivised adjoint orbit or in a bundle over such an orbit.  We then
have to classify which adjoint orbits in~$\mathfrak g^{\mathbb C}$ admit
$G$~actions of low cohomogeneity, before proceeding further.

The following constructions are essentially due to
Lichnerowicz~\cite{Lichnerowicz:complex}, but our treatment follows the
notation of~\cite{Swann:Trieste,Swann:HTwNil}.

The twistor space~$Z$ of a quaternionic K\"ahler manifold is a complex
contact manifold~\cite{Salamon:Invent}.  This means that $Z$~has a complex
line bundle~$L$ and a holomorphic one-form~$\theta\in\Omega^1(Z,L)$ such
that $\theta\wedge(d\theta)^n$ is nowhere zero.  The action of~$G^{\mathbb
C}$ preserves the complex contact structure.  We may define a \emph{moment
map} $f\colon Z\to L\otimes \mathfrak (g^{\mathbb C})^*$ for the action
of~$G^{\mathbb C}$ by
\begin{equation*}
  f(z)(Y)=\theta_z(\alpha(Y)), 
\end{equation*}
for $z\in Z$ and $Y\in\mathfrak g$.  Choosing an $\Ad$-invariant
inner-product on~$\mathfrak g^{\mathbb C}$, we identify $\mathfrak
g^{\mathbb C}$ with its dual equivariantly and regard $f$ as a map $Z\to
L\otimes\mathfrak g^{\mathbb C}$.

For~$z$ in~$U_0$ the vectors~$\alpha_z(Y)$ span~$T_zZ$, so non-degeneracy
of~$\theta$ implies that the map $f$~is non-vanishing at~$z$.
Projectivising, we get a well-defined map $\mathbb Pf\colon U_0\to \mathbb
P(\mathfrak g^{\mathbb C})$.  This map is $G^{\mathbb C}$-equivariant and
so the image~$\mathbb Pf(U_1)$ of the open orbit~$U_1$ lies in some
$G^{\mathbb C}$-orbit~$\mathbb P\mathcal O$.

The group~$\mathbb C^*$ acts on~$\mathfrak g^{\mathbb C}$ by scaling:
$X\mapsto \lambda X$, for $X\in\mathfrak g^{\mathbb C}$ and
$\lambda\in\mathbb C^*$.  We let $\mathbb P\colon \mathfrak g^{\mathbb
C}\setminus\{0\}\to\mathbb P(\mathfrak g^{\mathbb C})$ denote the
associated quotient map.  The adjoint action of~$G^{\mathbb C}$ commutes
with the scaling action of~$\mathbb C^*$, so the orbit~$\mathbb P\mathcal
O$ is $\mathbb P(\mathcal O)$ for some adjoint orbit~$\mathcal O$.

\begin{lemma}
  \label{lem:scaling}
  Let $\mathcal O$ be an adjoint orbit in~$\mathfrak g^{\mathbb
  C}\setminus\{0\}$.  If $\mathcal O$ is nilpotent, then $\mathcal O$~is
  invariant under the scaling action of~$\mathbb C^*$.  If $\mathcal O$~is
  not nilpotent, then the map $\mathcal O\to \mathbb P(\mathcal O)$ is an
  unbranched finite cover.
\end{lemma}

\begin{proof}
  Let $X$~be an element of~$\mathcal O$.  For $\lambda\in\mathbb C^*$, the
  eigenvalues of $\ad(\lambda X)$ are $\lambda$~times the eigenvalues
  of~$X$.  But if $\lambda X$ lies in~$\mathcal O$, the orbit of~$X$, then
  $\ad(\lambda X)$ has the same eigenvalues as~$\ad(X)$.  Thus either all
  the eigenvalues of~$\ad(X)$ are zero and $X$~is nilpotent, or
  $\abs\lambda=1$.  In the first case, it is well-known that nilpotent
  orbits are invariant under scaling: choose an $\Sl(2,\mathbb C)$-triple
  $\langle X,Y,H\rangle$ containing~$X$, then the action of~$\exp(tH)$
  on~$X$ gives the ray through~$X$ (see Carter~\cite{Carter:finite} for
  more details; recall that an $\Sl(2,\mathbb C)$-triple satisfies
  $[X,Y]=H$, $[H,X]=2X$ and $[H,Y]=-2Y$).  If $X$~is not nilpotent, then as
  $\ad(X)$~only has finitely many eigenvalues, there are only a finite
  number of possible $\lambda$ such that $\lambda X\in\mathcal O$.  Thus
  the map $\mathcal O\to \mathbb P(\mathcal O)$ is finite-to-one.  It is
  unbranched as it is $G^{\mathbb C}$-equivariant and the image is
  homogeneous.
\end{proof}

As $G$~is reductive, each adjoint orbit~$\mathcal O\subset\mathfrak
g^{\mathbb C}$ carries a complex symplectic structure~$\omega_{\mathcal
O}$.  This form was defined by Kirillov, Kostant and Souriau (see for
example~\cite{Kirillov:elements}) and is given by
\begin{equation*}
  \omega_{\mathcal O}([X,A],[X,B]) = \langle X, [A,B]\rangle,
\end{equation*}
where $\langle\cdot,\cdot\rangle$ is an invariant inner product.  In
particular, the complex dimension of each adjoint orbit~$\mathcal O$ is
even.  The previous lemma shows that the complex dimension of~$\mathbb
P(\mathcal O)$ is odd if $\mathcal O$ is nilpotent, and even otherwise.

Lichnerowicz shows that $\mathbb Pf$ can behave in two possible ways.  In
the proper case, $\mathcal O$ is nilpotent, so $\mathbb P(\mathcal O)$ is
contact. Moreover, $\mathbb Pf$ maps $U_1$ to $\mathbb P(\mathcal O)$ with
discrete fibre.

In the non-proper case, $\mathcal O$ is non-nilpotent, and $\mathbb Pf
(\mathcal O)$ is of even complex dimension. Now $\mathbb Pf$ has fibres of
complex dimension one with a complex Lie group structure.  Moreover, the
fibres are transverse to the complex contact distribution, so we have a
vector field $X$ on $U$ such that $\theta(X)$ never vanishes.  The bundle
$\mathcal L$ is therefore holomorphically trivial on $U_1$.
  
We shall now obtain estimates on the cohomogeneity of the adjoint
orbit~$\mathcal O$.

\begin{lemma}
  \label{lem:subman}
  Suppose the scalar curvature of~$(M,g)$ is positive. Then, in the
  non-proper case, $U_1$~contains no compact complex submanifolds.
\end{lemma}

\begin{proof}
  If we are in the non-proper situation, $\mathcal L$~is trivial on~$U_1$.
  However, if $(M,g)$~has positive scalar curvature, then $\mathcal L$~has
  a positive curvature form.  Hence $U_1$~contains no compact complex
  manifolds.
\end{proof}

In particular, the proper case must occur if $U_1$ contains a twistor line.

\begin{proposition}
  Suppose $(M,g)$ has positive scalar curvature and let $\mathcal
  O\subset\mathfrak g^{\mathbb C}$ be an orbit such that $\mathbb
  Pf(U_1)\subset\mathbb P(\mathcal O)$.
  
  If we are in the proper case, then $\mathcal O$~is a nilpotent orbit and
  has cohomogeneity at most~$5$ with respect to~$G$.  On the other hand, in
  the non-proper case, $\mathcal O$~is not nilpotent and has cohomogeneity
  at most~$2$ with respect to~$G$.
\end{proposition}

\begin{proof}
  Let us first consider the proper case.  As explained above, the map
  $\mathbb Pf \colon U_1 \mapsto \mathbb P(\mathcal O)$ is equivariant with
  discrete fibre.  Moreover the map $\mathbb P\colon \mathcal O \mapsto
  \mathbb P (\mathcal O)$ is equivariant with fibre $\mathbb C^*$.
  Therefore
  \begin{equation*}
    \cohom_G(Z)=\cohom_G U_1 = \cohom_G \mathbb P(\mathcal O) \geqslant
    \cohom_G(\mathcal O) -2
  \end{equation*}
  Now $\cohom_G(M)=1$, so $\cohom_G(Z) \leqslant 3$ and the above
  inequalities show that $\cohom_G(\mathcal O) \leqslant 5$, as required.
  
  For the non-proper case, using Lemma~\ref{lem:subman} we see that the
  fibres of~$\mathbb Pf\colon U_1\to\mathbb P(\mathcal O)$ are non-compact.
  We therefore have
  \begin{align*}
    \cohom_G (\mathcal O) &= \cohom_G(\mathbb P(\mathcal O))\\
    &\leqslant \cohom_G(U_1)-1 = \cohom_G(Z)-1 \leqslant 2.
  \end{align*}
\end{proof}

\begin{remark}
  The classification of these orbits is also relevant for the
  classification of $3$-Sasakian manifolds of cohomogeneity one and
  hyperK\"ahler manifolds of cohomogeneity two, see \S\S\ref{sec:hK-class}
  and~\ref{sec:3S-class}.
\end{remark}

\section{Cohomogeneity of Adjoint Orbits}
\label{sec:orbits}
In this section we shall first assume that $G$~is a compact simple Lie
group and give some procedures for calculating the cohomogeneity of adjoint
orbits in~$\mathfrak g^{\mathbb C}$ and related objects.  We will use this
to classify the orbits that arise from the discussion in the previous
section.  Results for semi-simple~$G$ will be given later
in~\S\ref{sec:ss-G}.

The adjoint orbits fall into three classes: semi-simple, nilpotent and
mixed.  We will discuss the cohomogeneity questions in each of these cases
in turn.  One feature all cases have in common is the existence of
monotonicity theorems for the cohomogeneity of topologically related
orbits.  In the nilpotent case, there is one quaternionic K\"ahler metric
known related to each orbit; we prove similar results for these structures.

In~\cite{Dancer-Swann:hK-cohom1}, the adjoint orbits of cohomogeneity one
were classified.  These were found to be the minimal nilpotent orbits and
the semi-simple orbit of $\diag(\lambda,\dots,\lambda,-n\lambda)$
in~$\Sl(n+1,\mathbb C)$; this last orbit is $SU(n+1)$-equivariantly
diffeomorphic to $T^*\CP(n)$.

\subsection{Semi-Simple Orbits}
\label{sec:ss}
Let $\mathcal O$~be a semi-simple orbit.  In this case, $\mathcal O$~is a
flag manifold and we have $G$-equivariant diffeomorphisms
\begin{equation*}
  \mathcal O
  \cong T(G/K)
  \cong G^{\mathbb C}/K^{\mathbb C},
\end{equation*}
where $K$~is the centraliser of some torus in~$G$.  We will write
$\mathfrak g\cong\mathfrak k\oplus\mathfrak m$ with $\mathfrak m$ a
$K$-module.

First we note

\begin{lemma}
  \label{lem:m}
  If $G/K$~is any compact homogeneous space, then
  \begin{equation*}
    \cohom_G T(G/K) = \cohom_K \mathfrak m
  \end{equation*}
  and the latter is bounded below by the number of irreducible summands
  in~$\mathfrak m$. \qed
\end{lemma}

We will obtain bounds on the cohomogeneity of flag manifolds by considering
various fibrations.  The following will be a useful observation.

\begin{lemma}
  \label{lem:fibration}
  Suppose $G/K \to G/H$ is a fibration of homogeneous $G$-spaces, with $G$,
  $K$ and~$H$ compact and $\dim K < \dim H$.  Then
  \begin{equation*}
    \cohom T(G/K) > \cohom T(G/H).
  \end{equation*}
\end{lemma}

\begin{proof}
  We have $\mathfrak g \cong \mathfrak h\oplus \mathfrak n$ with $\mathfrak
  k\subsetneqq\mathfrak h$.  So $\mathfrak m \cong \mathfrak a\oplus
  \mathfrak n$ for some non-zero representation~$\mathfrak a$ of~$K$.  As
  $K$~is compact, this gives
  \begin{equation*}
    \cohom_K \mathfrak m > \cohom_K \mathfrak n \geqslant \cohom_H
    \mathfrak n,
  \end{equation*}
  since each $K$-orbit in~$\mathfrak n$ lies in an $H$-orbit.  The result
  now follows from Lemma~\ref{lem:m}.
\end{proof}

Recall that a flag manifold~$T(G/K)$ is specified by a subset~$\mathfrak K$
of the simple roots for~$\mathfrak g$: the elements of~$\mathfrak K$ are
the simple roots not in~$\mathfrak k^{\mathbb C}$.  A Dynkin diagram is
defined for~$G/K$ by taking the Dynkin diagram for~$\mathfrak g$ and
putting a cross through each simple root~$\alpha$ in~$\mathfrak K$
(see~\cite{Baston-Eastwood:Penrose}, for example).  We call $\abs{\mathfrak
K}$ the \emph{length} of~$G/K$ and note that there is a
fibration~$T(G/K_1)\to T(G/K_2)$ if and only if $\mathfrak K_1\supset
\mathfrak K_2$.  The following result is a simple consequence of the
previous Lemma.

\begin{lemma}[Monotonicity for Semi-Simple Orbits]
  \label{lem:ss-mono}
  If $\mathfrak K_1\supset\mathfrak K_2$, then
  \begin{equation*}
    \cohom T(G/K_1) \geqslant \cohom T(G/K_2) +
    \left(
      \abs{\mathfrak K_1} - \abs{\mathfrak K_2}
    \right).
  \end{equation*}
  In particular, $\cohom T(G/K_1) \geqslant \abs{\mathfrak K_1}$. \qed
\end{lemma}

This last estimate is extremely crude.  The next result, which is
particularly relevant because of our interest in flag manifolds of
cohomogeneity two, serves to emphasise this.

\begin{lemma}
  \label{lem:length2}
  If $T(G/K)$ is a flag manifold with $\abs{\mathfrak K}\geqslant 2$, then
  \begin{equation*}
    \cohom T(G/K) \geqslant 3.
  \end{equation*}
\end{lemma}

\begin{proof}
  Suppose $\alpha_1$ and $\alpha_2$ are simple roots in~$\mathfrak K$.  Let
  $\alpha$~denote the top root form of~$\mathfrak g$.  Then $\alpha_1$,
  $\alpha_2$ and $\alpha$~do not pairwise differ by linear combinations of
  roots of~$\mathfrak k^{\mathbb C}$.  Therefore, by Kostant's
  criterion~\cite[p.~40]{Black:harmonic}, $\alpha_1$, $\alpha_2$ and
  $\alpha$ lie in distinct $K$-submodules of~$\mathfrak m$ (using the
  notation of Lemma~\ref{lem:m}).  In particular, $\mathfrak m$~has at least
  three summands and hence $\cohom_K\mathfrak m\geqslant 3$, giving the
  result.
\end{proof}

The above results now enable us to classify semi-simple orbits of
cohomogeneity two.

\begin{theorem}
  \label{thm:ss-c2}
  For a compact simple group~$G$, the semi-simple orbits of cohomogeneity
  two are the tangent bundles of the following homogeneous spaces:
  \begin{gather*}
    \frac{\SP(n+1)}{\Un(1)\SP(n)},\qquad
    \frac{\SU(n+2)}{\Special(\Un(n)\times\Un(2))}\\
    \frac{\SO(n+2)}{\SO(n)\times\SO(2)}, \qquad
    \frac{\SO(10)}{\Un(5)}
    \qquad\text{and}\qquad
    \frac{\LieE_6}{\Spin(10)\SO(2)}.
  \end{gather*}
\end{theorem}

\begin{proof}
  If $T(G/K)$~is a flag manifold, Burstall \&\ Rawnsley
  \cite{Burstall-Rawnsley:twistor} show that there is a fibration $G/K\to
  G/H$ with $G/H$~an inner symmetric space.  We can take $H$ to be
  connected.  By Lemma~\ref{lem:length2}, $T(G/K)$~is a flag manifold of
  length one in our case.  So Lemma~\ref{lem:fibration} implies we have two
  cases to consider, either the fibration~$G/K\to G/H$ is non-trivial or
  $G/K$~is Hermitian symmetric.
  
  If the fibration~$G/K\to G/H$ is non-trivial, then $T(G/H)$~must have
  cohomogeneity one.  This implies that $G/H$~is rank
  one (see~\cite{Besse:Einstein} or the remarks in the last paragraph of the
  proof).
  Moreover, $G/H$ cannot be Hermitian
  symmetric, otherwise the length of~$T(G/K)$ would be at least two.  It
  follows that $G/H$~is one of
  \begin{gather*}
    \text{(i)}\quad \frac{\SP(n+1)}{\SP(1)\SP(n)},\qquad \text{(ii)}\quad
    \frac{\SO(2n+1)}{\SO(2n)},\qquad \text{(iii)}\quad
    \frac{\LieF_4}{\Spin(9)}.
  \end{gather*}
  The rank one symmetric space~$\SO(2n)/\SO(2n-1)$ does not occur, as it is
  not inner.
  
  In case~(i), there are two flag manifolds of length one fibering
  over $G/H$, namely $\Dyx\Dyo\cdots\Dyo\ltwoa\Dyo$ and
  $\Dyo\Dyx\Dyo\cdots\Dyo\ltwoa\Dyo$.  For the first of these
  $K=\Un(1)\SP(n)$ and $\mathfrak m^{\mathbb C}\cong (L+\overline L)E + L^2
  + {\overline L}^2$, where $E\cong \mathbb C^{2n}$ and $L\cong\mathbb C$
  are the standard representations of~$\SP(n)$ and $\Un(1)$~respectively.
  In particular, $\mathfrak m\cong\mathbb H^n\oplus\mathbb R^2$, with
  $\SP(n)$ acting transitively on the unit sphere in~$\mathbb H^n$ and
  $\Un(1)$~acting non-trivially on~$\mathbb R^2$.  So this is of
  cohomogeneity two and we have the first case on our list.

  For~$\Dyo\Dyx\Dyo\cdots\Dyo\ltwoa\Dyo$, we have $K=\SU(2)\Un(1)\SP(n-1)$
  and $\mathfrak m^{\mathbb C}\cong (L^2+{\overline L}^2)S^2H+(L+\overline
  L)HE$, with $H\cong\mathbb C^2$ the standard representation of~$\SU(2)$
  and $E$~and $L$ essentially as before.  The real module underlying
  $(L^2+{\overline L}^2)S^2H$ is~$\mathbb R^6$ with a cohomogeneity two
  action of~$\SO(2)\times\SO(3)$.  Thus, $\mathfrak m$~has cohomogeneity at
  least three under~$K$ and so is not on our list.
  
  For case~(ii), we take $n>2$ to avoid overlapping with case~(i).  Then
  the only flag manifold of length one fibering over this is
  $\Dyo\cdots\Dyo\rtwoa\Dyx$, which has $K=\Un(n)$ and $\mathfrak
  m\cong[\Lambda^{1,0}\mathbb C^n]+[\Lambda^{2,0}\mathbb C^n]$, and which
  is not of cohomogeneity two.  (The brackets~$[\cdot]$ indicate the
  underlying real representation.)
  
  Finally, for case~(iii), first note that $\spin(9)$ is embedded
  in~$\lieF_4$ with positive simple roots $01\rtwoa22$, $10\rtwoa00$,
  $01\rtwoa00$, $00\rtwoa10$.  The only flag manifold of length one
  fibering over $\LieF_4/\Spin(9)$ is~$\Dyo\Dyo\rtwoa\Dyo\Dyx$.  Here
  $\mathfrak k\cong\un(1)\oplus\so(7)$, which is $22$-dimensional, whereas
  $\dim\mathfrak m=30$, so the action is not cohomogeneity two.
  
  We now need to consider the case when $G/K$~is a Hermitian symmetric
  space.  The cohomogeneity we are interested in is $\cohom_K \mathfrak m$,
  where $\mathfrak m$ is the isotropy representation for this space.  The
  Cartan theory for symmetric spaces (see \cite[Chapter
  VIII]{Helgason:symmetric}, for example) shows that a transversal for the
  action of~$K$ on~$\mathfrak m$ is given by a chamber in a maximal Abelian
  subspace of~$\mathfrak m$.  It follows that~$\cohom_K \mathfrak m$ is
  just the rank of the symmetric space~$G/K$.  The table in Besse
  \cite[pp.\ 312--313]{Besse:Einstein} shows that the Hermitian symmetric
  spaces of rank two are (modulo low-dimensional coincidences, and up to
  covers) the following: the Grassmannians of two-planes in~$\mathbb C^n$,
  the hyperquadrics $\SO(n+2)/\SO(n) \times \SO(2)$, the symmetric space
  $\SO(10)/\Un(5)$ and the exceptional space $E_6/\Spin(10) \SO(2)$.  This
  completes the proof.
\end{proof}

For future reference, it is worth noting that the principal orbits
in~$\mathcal O$ for the five cases of Theorem~\ref{thm:ss-c2} are
\begin{gather}
  \frac{\SP(n+1)}{\SP(n-1)} , \qquad \frac{\SU(n+2)}{\Special(\Un(n-2)
  \times  \Un(1) \times \Un(1))},\notag \\
  \frac{\SO(n+2)}{\SO(n-2)}, \qquad \frac{\SO(10)}{\Special(\Un(2) \times
  \Un(2))} \qquad\text{and}\qquad \frac{\LieE_6}{\Un(4)}.
  \label{eq:ss-p-orbit-1}
\end{gather}
Also, from the cohomogeneity one case~$T^*\CP(n)$ we have the principal
orbit 
\begin{equation}
  \label{eq:ss-p-orbit-2}
  \frac{\SU(n+1)}{\Un(n-1)}
\end{equation}
However, not all these cases will arise from the twistor spaces of
cohomogeneity one quaternionic K\"ahler manifolds, as we will see
in~\S\ref{sec:qK-class}.

\subsection{Nilpotent Orbits}
\label{sec:nil}
Recall that there is a partial order on the nilpotent orbits in~$\mathfrak
g^{\mathbb C}$ defined by $\mathcal O_1 \succeq \mathcal O_2$ if and only
if $\overline{\mathcal O_1} \supset \mathcal O_2$.

\begin{proposition}[Monotonicity for Nilpotent Orbits]
  \label{prop:nil-mon}
  If $\mathcal O_1$ and $\mathcal O_2$ are nilpotent orbits with $\mathcal
  O_1 \succneqq \mathcal O_2$, then
  \begin{equation*}
    \cohom \mathcal O_1 > \cohom \mathcal O_2.
  \end{equation*}
\end{proposition}

\begin{proof}
  Fix an element~$X$ in~$\mathcal O_2$.  For any $\Sl(2,\mathbb C)$-triple
  $\langle X, Y, H\rangle$, set
  \begin{equation*}
    S_{X,Y} = X + \mathfrak z(Y),
  \end{equation*}
  where $\mathfrak z(Y)$~denotes the centraliser of~$Y$.
  Slodowy~\cite{Slodowy:singularities} showed that $S_{X,Y}$~is a transverse
  slice to~$\mathcal O_2$ at~$X$.  We need to study the intersection of
  $G$-orbits with such a slice.  For this we make a good choice of~$Y$.

  Let $B=\stab_G X$ be the real stabiliser of~$X$ and let $\mathfrak c$~be
  the centraliser in~$\mathfrak g^{\mathbb C}$ of~$\mathfrak b$, the Lie
  algebra of~$B$.  Then $\mathfrak c$~is $\sigma$-invariant and hence
  reductive.  Moreover, $\mathfrak c$ contains~$X$, so there is a
  $Y\in\mathfrak c$ such that $\langle X,Y,[X,Y]\rangle$ is an
  $\Sl(2,\mathbb C)$-triple.  As $Y$~lies in~$\mathfrak c$, we have that
  $Y$, $\mathfrak z(Y)$ and $S_{X,Y}$ are all $B$-invariant.
  
  Fix this choice of~$Y$ and let $H=[X,Y]$.  In $\mathfrak c$, $\langle
  X,Y,H\rangle$~is conjugate to a real $\Sl(2,\mathbb C)$-triple, so
  there exists $g\in \exp\mathfrak c$ such that $(\Ad g)Y=-\sigma((\Ad
  g)X)$.  Thus for $\sigma^g=(\Ad g)^{-1}\sigma(\Ad g)$, we have
  $Y=-\sigma^g X$ and $\sigma^g\mathfrak b=\mathfrak b$.  Define $G^g$~to
  be the connected subgroup of~$G^{\mathbb C}$ with Lie algebra $\mathfrak
  g^g = \mathfrak g^{\mathbb C} \cap \sigma^g \mathfrak g^{\mathbb C}$.
  Then $G^g$~is isomorphic to~$G$ and contains~$B$.

  Choose a $\sigma^g$-invariant $\Sl(2,\mathbb C)$-triple $\langle
  X_1,Y_1,H_1\rangle$ with $X_1\in \mathcal O_1$ and $Y_1=-\sigma^g X_1$.
  Let $\mathcal M$ be the space of maps $(X(t),H(t))\colon \mathbb R \to
  \mathfrak g^{\mathbb C}\times i\mathfrak g^g$ such that
  \begin{enumerate}
  \item[(a)] $\dot X=-2X+[X,H]$, $\dot H=-2H+2[X,\sigma^gX]$,
  \item[(b)] $X(t)\to X$ and $H(t)\to H$ as $t\to+\infty$, and
  \item[(c)] $X(t)\to X_{-\infty}$ and $H(t)\to H_{-\infty}$ as
    $t\to-\infty$, with $(X_{-\infty},H_{-\infty})$ $G^g$-conjugate
    to~$(X_1,H_1)$.
  \end{enumerate}
  Kronheimer~\cite{Kronheimer:nilpotent} showed that $\mathcal M$~is
  a manifold naturally isomorphic to~$S_{X,Y} \cap \mathcal O_1$.  This
  isomorphism is $\stab_{G^g} X$ equivariant; but $B\subset \stab_{G^g} X$,
  so we have $\mathcal M\cong S_{X,Y} \cap \mathcal O_1$
  $B$-equivariantly.
  
  The manifold~$\mathcal M$ admits an action of~$\mathbb R$ given by
  $(X(t), H(t)) \mapsto (X(t-c), H(t-c))$.  This $\mathbb R$-action is
  non-trivial and commutes with the action of~$G^g$ on the values.  In
  particular, $\mathcal M/B$~has positive dimension.

  Let $\mathfrak v$~be a $B$-invariant complement to~$T((\Ad G)X)$
  in~$T_X\mathcal O_2$.  Then the orbit~$(\Ad G)X$ has a $G$-invariant
  neighbourhood
  \begin{equation*}
    U \cong G \times_B (\mathfrak v+\mathfrak z(Y))
  \end{equation*}
  in~$\mathfrak g^{\mathbb C}$ such that the map $U\cap\overline {\mathcal
  O_1}\to G\times_B\mathfrak v$ is an equivariant surjection onto a
  $G$-invariant neighbourhood of~$(\Ad G)X$ in~$\mathcal O_2$ and contains
  $S_{X,Y}\cap \mathcal O_1$ in the fibre over~$0\in\mathfrak v$.  Thus
  for generic~$X$ in~$\mathcal O_2$, we have
  \begin{equation*}
    \begin{split}
      \cohom_G \mathcal O_1
      &= (\cohom_B S_{X,Y}\cap \mathcal O_1) + \cohom_G \mathcal O_2 \\
      &> \cohom_G \mathcal O_2,
    \end{split}
  \end{equation*}
  as required.
\end{proof}

Each simple Lie group~$G$ has a unique non-trivial nilpotent orbit of
smallest dimension.  We call this the \emph{minimal} nilpotent
orbit~$\mathcal O_{\text{min}}$.  The orbit $\mathcal O_{\text{min}}$ has
the property that $\mathcal O_{\text{min}} \preceq \mathcal O$ for all
non-trivial nilpotent orbits~$\mathcal O$.
In~\cite{Dancer-Swann:hK-cohom1}, it was shown that $\mathcal
O_{\text{min}}$~is the unique nilpotent orbit of cohomogeneity one.

As explained in~\S\ref{sec:moment}, each nilpotent orbit~$\mathcal O$ fibres
over a quaternionic K\"ahler manifold~$\mathfrak M(\mathcal O)=\mathcal
O/\mathbb H^*$.  In the case of~$\mathcal O_{\text{min}}$, the space
$\mathfrak M(\mathcal O_{\text{min}})$~is $G$-homogeneous and is the Wolf
space with isometry group~$G$.  (The details of the fibration in this case
may be found in~\cite[\S6]{Swann:MathAnn}.)  Proposition~\ref{prop:nil-mon}
carries over to the manifolds~$\mathfrak M(\mathcal O)$.

\begin{proposition}[Monotonicity for $\mathfrak M(\mathcal O)$]
  \label{prop:nqK-mon}
  If $\mathcal O_1$ and $\mathcal O_2$ are nilpotent orbits with $\mathcal
  O_1 \succneqq \mathcal O_2$, then
  \begin{equation*}
    \cohom \mathfrak M(\mathcal O_1) > \cohom \mathfrak M(\mathcal O_2).
  \end{equation*}
\end{proposition}

\begin{proof}
  Fix $X \in \mathcal O_2$ and choose $Y$~as in the proof of
  Proposition~\ref{prop:nil-mon}.  The subgroup $\mathbb R^* \leqslant
  \mathbb H^*$ acts on~$X$ by $X \mapsto \lambda^2 X$ and an $\Sl(2,\mathbb
  C)$-triple containing~$\lambda^2 X$ is given by $\langle \lambda^2 X,
  \lambda^{-2} Y, [X,Y]\rangle$.  Now $\mathfrak z(\lambda^{-2} Y) =
  \mathfrak z(Y)$, so the slice at~$\lambda^2 X$ is given by~$S_{\lambda^2
  X,Y}$.  Note that the set $S_{X,Y} \cap S_{\lambda^2 X,Y}$ is empty
  for~$\lambda^2 \ne 1$.  We set 
  \begin{equation*}
    W_{X,Y} = \mathbb H^*(S_{X,Y} \cap \overline{\mathcal O_1} ) / \mathbb
    H^* = \SP(1) (S_{X,Y} \cap \overline{\mathcal O_1} ) / \SP(1)
  \end{equation*}
  to get a slice to $\mathfrak M(\mathcal O_2) \subset \overline{\mathfrak
  M(\mathcal O_1)}$ at~$\mathbb H^*X$.

  The proof of Proposition~\ref{prop:nil-mon} shows that $S_{X,Y} \cap
  \mathcal O_1$ is non-compact.  As $\SP(1)$ and the stabiliser~$D$
  of~$\mathbb H^*X$ are both compact, we have for generic $X \in \mathcal
  O_2$,
  \begin{equation*}
    \begin{split}
      \cohom_G \mathfrak M(\mathcal O_1)
      &= \cohom_D W_{X,Y} + \cohom_G \mathfrak M(\mathcal O_2) \\
      &> \cohom_G \mathfrak M(\mathcal O_2),
    \end{split}
  \end{equation*}
  as required.
\end{proof}

{\tolerance=250
Let us now consider the cohomogeneity-one quaternionic K\"ahler manifolds
that arise from nilpotent orbits.  We say that $\mathcal O$ is next to
minimal in the partial order~$\preceq$ if $\mathcal O\succneqq \mathcal
O_{\text{min}}$ and there are no nilpotent orbits $\mathcal O'$ with
$\mathcal O \succneqq \mathcal O' \succneqq \mathcal O_{\text{min}}$.  The
Lie algebra~$\mathfrak g^{\mathbb C}$ may have more than one nilpotent
orbit satisfying this condition, as we will see in the proof of the
following result.
\par
}

\begin{theorem}
  \label{thm:nqK-c1}
  Suppose $\mathcal O$ is a nilpotent orbit and $\mathfrak M(\mathcal
  O)$~is the associated quaternionic K\"ahler manifold.  Then $\mathfrak
  M(\mathcal O)$~is of cohomogeneity one if and only if $\mathcal O$ is
  next to minimal in the order~$\preceq$ on nilpotent orbits.
\end{theorem}

\begin{proof}
  Proposition~\ref{prop:nqK-mon} implies that the only candidates are the
  next to minimal nilpotent orbits.  We will verify that each of these does
  indeed lead to a cohomogeneity-one structure on~$\mathfrak M(\mathcal O)$
  by considering each case in turn.
  
  First, we recall from~\cite{Swann:HTwNil} the structure of~$\mathfrak
  M(\mathcal O)$ as a $G$-manifold.  The nilpotent orbit~$\mathcal O$
  defines a conjugacy class of subalgebras~$\Sl(2,\mathbb C) \leqslant
  \mathfrak g^{\mathbb C}$.  Choose one such subalgebra~$\Sl(2,\mathbb C)$
  and without loss of generality, assume that this is the complexification
  of a real subalgebra~$\su(2)$.  Let $K$~be the centraliser
  of~$L=\exp\{\su(2)\}$ in~$G$.  We may decompose $\mathfrak g$ under the
  adjoint action of~$L\times K$ to get
  \begin{equation*}
    \mathfrak g \cong \su(2) \oplus \mathfrak
    k \oplus \bigoplus_{k\geqslant 1} [A_kS^k],
  \end{equation*}
  where $S^k$~is the irreducible representation of~$\SU(2)$ of complex
  dimension~$k+1$ and $A_k$~is a $K$-module.  Then $\mathfrak M(\mathcal
  O)$ is $G$-equivariantly isomorphic to the bundle 
  \begin{equation*}
    W=\bigoplus_{k\geqslant 2} [A_kS^{k-2}] \longrightarrow G/LK.
  \end{equation*}
  In particular, when $\mathfrak M(\mathcal O)$~is of cohomogeneity one,
  there is at most one non-zero module~$A_k$ for $k\geqslant 2$.
  
  We now consider the various possible types of~$G$ in turn.  Relevant
  facts about the description of nilpotent elements and their stabilisers
  have been conveniently collected in~\cite[Chapter~13]{Carter:finite}.
  
  Suppose $G$~is of type~$A_n$.  Then the nilpotent orbits are described by
  partitions of~$n+1$ giving the sizes of the Jordan blocks.  The minimal
  nilpotent orbit is~$(21^{n-1})$ and there is a unique next-to-minimal
  nilpotent orbit.  When $n>2$, the next-to-minimal orbit is~$(2^21^{n-3})$
  and for this case the Lie algebra of~$K$ is~$\su(2)_-+\su(n-3)+\un(1)$.
  Write the Lie algebra of~$L$ as~$\su(2)_+$ to distinguish it from the
  $\su(2)$-factor in~$\mathfrak k$.  Now as a $KL$-module, $\mathbb
  C^{n+1}\cong S_+S_-\ell^r+V\ell^s$, where $S_\pm\cong\mathbb C^2$,
  $V\cong\mathbb C^{n-3}$ and~$\ell\cong\mathbb C$.  So
  \begin{equation*}
    \begin{split}
      \su(n+1)
      &= \End_0\mathbb C^{n+1} \\
      &\cong \su(2)_+ + \mathfrak k + [S_+S_-\ell^{r-s}V^*] + [S^2_+S^2_-].
    \end{split}
  \end{equation*}
  Thus, $W=[S^2_-]$ which is~$\mathbb R^3$ with the standard
  representation of~$\SO(3)$.  Thus $\mathfrak M(2^21^{n-3})$~is
  cohomogeneity one.
  
  For $n=2$, there are only two non-trivial nilpotent orbits.  The
  next-to-minimal orbit is the regular orbit.  We have
  $\su(3)=\su(2)+[S^4]$ and $W=[S^2]$, which is cohomogeneity one.
  
  For type~$B_n$, the nilpotent orbits are again described by partitions,
  this time of~$2n+1$, but not all partitions arise.  The minimal orbit
  is~$(2^21^{2n-3})$.  There are two next-to-minimal orbits: $(31^{2n-2})$
  and $(2^41^{2n-7})$.  For the first of these, $\mathfrak k\cong
  \so(2n-2)$.  The $\SO(2n+1)$-module~$\mathbb R^{2n+1}$ splits as
  $[S^2]+V$, with $V=\mathbb R^{2n-2}$ the standard representation of~$K$.
  Using $\so(2n+1,\mathbb C)\cong\Lambda^2\mathbb C^{2n-2}$, we find
  $W\cong V$, which is cohomogeneity one.  For~$(2^41^{2n-7})$, $\mathfrak
  k\cong\LSP(2)+\so(2n-7)$ and $\mathbb R^{2n+1}\cong [\mathbb
  C^4S^1]+\mathbb R^{2n-7}$, with $\SP(2)$ and $\SO(2n-7)$~acting
  irreducibly on~$\mathbb C^4$ and~$\mathbb R^{2n-7}$, respectively.  This
  gives $W\cong\mathbb R^5$ with $\SP(2)$-acting via the standard
  representation of~$\SO(5)$, which is again cohomogeneity one.
  
  The nilpotent orbits for $G$ of type~$C_n$ are classified by certain
  partitions of~$2n$.  The minimal orbit is~$(21^{2n-2})$.  There is only
  one next to minimal orbit: $(2^21^{2n-4})$.  The centraliser~$K$ has Lie
  algebra~$\LSP(n-2)+\un(1)$ and $\mathbb C^{2n}\cong
  (\ell+\ell^{-1})S^1+\mathbb C^{2n-4}$, where $\LSP(n-2)$~acts irreducibly
  on the last summand and $\ell$~is a one-dimensional representation
  of~$\Un(1)$.  We have $\LSP(n,\mathbb C)\cong S^2\mathbb C^{2n}$ and
  hence $W\cong\mathbb R^2$ as an irreducible $\SO(2)$-module.  This is
  cohomogeneity one.
  
  For type~$D_n$, we have $G=\SO(2n)$, the minimal orbit is~$(2^21^{2n-4})$
  and the next to minimal orbits have partitions $(31^{2n-3})$ and
  $(2^41^{2n-8})$.  For $n$~odd this gives two next-to-minimal orbits,
  however for $n$~even, the last partition describes two orbits, which we
  need not distinguish, and there are three next-to-minimal orbits.  The
  calculations for all of these partitions are the same as the $B_n$~case
  and give cohomogeneity-one manifolds each time.
  
  For the exceptional groups there is a unique next-to-minimal orbit.  The
  centralisers, Lie algebra decompositions and representation~$W$ are given
  in Table~\ref{tab:ntme}.  The calculations for the groups of type~$E$
  where done using the program ``Lie''~\cite{Lie:v2}.  The script for the
  case of~$\LieE_8$ is given in Figure~\ref{fig:e8}.  This script first
  calculates the semi-simple element~\texttt h associated to the nilpotent
  orbit.  It then determines the stabiliser~$\texttt{cs}=\mathfrak k_1$
  of~\texttt h and finds its Cartan type.  The decomposition of the adjoint
  representation of~$\LieE_8$ under~$\mathfrak k_1$ is calculated, using
  the restriction matrix~\texttt{rm}.  This information, combined with
  direct calculation of the weight spaces for the action of~$\su(2)$, are
  enough to determine the where the nilpotent element~$X$ lies.  The Lie
  algebra of~$\mathfrak k$ is now the stabiliser in~$\mathfrak k_1$ of~$X$
  and decomposing each $\mathfrak k_1$-module under the action
  of~$\mathfrak k$ leads to~$W$.
\end{proof}

\begin{table}[htbp]
  \begin{center}
    \leavevmode
    \begin{tabular}{lcclc}
      \hline
      Type & Orbit & $\mathfrak k$ &  $W$ & $\cohom(\mathcal O)$\\
      \hline
      \hline
      \vrule height 12pt width 0pt $A_n$ \\
      $n=2$ & $(3)$ & $\{0\}$ & $[S^2] \cong \mathbb R^3$ & $4$ \\
      $n\geqslant3$ & $(2^21^{n-3})$ & $\su(2)_-+\su(n-3)$ &
      $[S^2_-] \cong \mathbb R^3$ & $2$
      \\
      & & $+ \un(1)$ \\
      \hline
      \vrule height 12pt width 0pt \rlap{$B_{(n-1)/2}$, $D_{n/2}$}\\
      & $(31^{n-3})$ & $\so(n-3)$ & $\mathbb R^{n-3}$ & $2$
      \\
      & $(2^41^{n-8})$ & $\so(5) + \so(n-8)$ & $\mathbb R^5$ & $2$ \\
      \hline
      \vrule height 12pt width 0pt $C_n$ & $(2^21^{2n-4})$ & $\so(2) +
      \LSP(n-2)$ & $\mathbb R^2$ & $2$ \\ 
      \hline
      \vrule height 12pt width 0pt $\LieG_2$ & ${\scriptstyle 0\rthreea1}$
      & $\su(2)_-$ 
      & $[S^1_+S^1_-]\cong\mathbb R^4$ & $2$ \\
      $\LieF_4$ & ${\scriptstyle 00\rtwoa01}$ &
      $\so(6)$ & $\mathbb R^6$ & $2$ \\
      $\LieE_6$ & $\scriptstyle 10\Dyatop0001$ & $\so(2)+\so(7)$ & $\mathbb
      R^7$ & $2$ \\ 
      $\LieE_7$ & $\scriptstyle 010\Dyatop0000$ & $\su(2)+\so(9)$ &
      $\mathbb R^9$ & $2$ \\
      $\LieE_8$ & $\scriptstyle 0000\Dyatop0001$ & $\so(13)$ &
      $\mathbb R^{13}$ & $2$ \\
      \hline
      \\
    \end{tabular}
    \caption{Next-to-minimal orbits, specified by partitions for classical
    groups and Dynkin diagrams for exceptional Lie groups.  In each case
    $W=\mathbb R^n$~is the standard irreducible representation
    of~$\SO(n)$.}
    \label{tab:ntme}
  \end{center}
\end{table}

\begin{figure}[htbp]
  \begin{center}
    \begin{minipage}{52.5ex}
      \leavevmode
      \begin{small}
\begin{verbatim}
setdefault E8
print(diagram)
h=i_Cartan*[1,0,0,0,0,0,0,0]/det_Cartan;print(h)
hp=h+0;print(hp)
s=cent_roots(hp)
ct=Cartan_type(s);print(ct)
cs=closure(s)
rm=res_mat(cs)
ad=expon(adjoint,1);print(ad);print(dim(ad))
branch(ad,ct,rm)
\end{verbatim}
      \end{small}
    \end{minipage}
    \caption{``Lie'' script for calculations for~$\LieE_8$ in the
    nilpotent case.} 
    \label{fig:e8}
  \end{center}
\end{figure}

Note that if we calculate the cohomogeneities of the next-to-minimal
nilpotent orbits themselves we do not quite get this uniform answer.
Having seen the above theorem one might guess that moving up the partial
order decreases the cohomogeneity of~$\mathfrak M(\mathcal O)$ by one at
each stage.  However, a look at the diagrams
in~\cite{Kraft-Procesi:classical} for~$\LSP(4)$, $\so(10)$ or $\su(7)$
shows that this can not be the case, as there can be paths of different
lengths in the partial order.

\subsection{Mixed Orbits}
\label{sec:mixed}
If $X$~is a general element of~$\mathfrak g^{\mathbb C}$, we can write $X$
uniquely as $X=X_s+X_n$ with $X_s$~semi-simple, $X_n$~nilpotent and
$[X_s,X_n]=0$.  Similarly, if $\mathcal O$ is the $G^{\mathbb C}$-orbit
of~$X$, then we write $\mathcal O_s$ and $\mathcal O_n$ for the orbits of
$X_s$ and~$X_n$, respectively.  We partially order the orbits by saying
$\mathcal O_1 \succeq \mathcal O_2$ if and only if $\overline{\mathcal O_1}
\supset \mathcal O_2$.  

\begin{lemma}
  \textup{(a)} If two mixed orbits $\mathcal O_1$ and~$\mathcal O_2$
  satisfy $\mathcal O_1 \succeq \mathcal O_2$ then $(\mathcal O_1)_s =
  (\mathcal O_2)_s$ and $(\mathcal O_1)_n \succeq (\mathcal O_2)_n$.
  
  \textup{(b)} If $\mathcal O$~is a mixed orbit, then $\mathcal O \succeq
  \mathcal O_s$.
\end{lemma}

\begin{proof}
  (a)
  Suppose $\mathcal O_2$~is in the closure of~$\mathcal O_1$.
  Fix $X\in\mathcal O_2$ and choose a sequence $X^{(i)}\in\mathcal O_1$
  such that $X^{(i)}\to X$.  Now the orbit of the semi-simple part
  of~$X^{(i)}$ is determined by the characteristic polynomial
  of~$\ad(X^{(i)})$.  But that polynomial is a continuous function
  of~$X^{(i)}$ and is constant on~$\mathcal O_1$.  So the characteristic
  polynomial of~$\ad(X)$ equals that of~$\ad(X^{(i)})$ and hence $(\mathcal
  O_1)_s=(\mathcal O_2)_s$.  The map $\pi_n\colon \overline{\mathcal O_1}
  \to \overline{(\mathcal O_1)_n}$ given by $\pi_n(X)=X_n=X-X_s$ is
  continuous and $G^{\mathbb C}$-equivariant.  So $\pi_n(\overline O_1)$ is
  contained in $\overline{\pi_n(\mathcal O_1)}=\overline{(\mathcal
  O_1)_n}$.  In particular, $(\mathcal O_2)_n$ is contained
  in~$\overline{(\mathcal O_1)_n}$, as required.
  
  (b) Write $X=X_s+X_n$.  Then $X_s+\lambda X_n$~is in~$\mathcal O$ for all
  non-zero $\lambda\in\mathbb C$.  Thus $X_s$~is in the closure
  of~$\mathcal O$.
\end{proof}

\begin{proposition}[Monotonicity for Mixed Orbits]
  \label{prop:mixed-mon}
  Suppose $\mathcal O_1$ and $\mathcal O_2$ are mixed orbits with $\mathcal
  O_1 \succneqq \mathcal O_2$.  Then
  \begin{equation*}
    \cohom \mathcal O_1 > \cohom \mathcal O_2.
  \end{equation*}
\end{proposition}

\begin{proof}
  This is a slight modification of the proof of
  Proposition~\ref{prop:nil-mon}.  Fix $X \in \mathcal O_2$.  The
  semi-simple orbit~$(\mathcal O_2)_s$ is isomorphic to~$T(G/K)$ for some
  $K \leqslant G$.  Suppose $X_s$~lies in the fibre of $T(G/K)\to G/K$ over
  $gK$.  Then $\stab_G X_s$~is a subgroup of~$gKg^{-1}$ which is a compact
  real form of~$g K^{\mathbb C} g^{-1}$.  Now $g K^{\mathbb C} g^{-1}$ is
  abstractly isomorphic to~$K_s^{\mathbb C} := \stab_{G^{\mathbb C}}
  X_s$, so we conclude that $\mathfrak k_s^{\mathbb C} = \mathfrak z(X_s)$
  admits a real structure~$\sigma_s$ giving a compact real form~$K_s$
  containing~$\stab_G X_s$.
  
  Let $B=\stab_G X=\stab_G X_s\cap\stab_G X_n$ and let $\mathfrak c$ be the
  centraliser of~$\mathfrak b$ in~$\mathfrak k_s^{\mathbb C}$.  Then $X$,
  $X_s$ and $X_n$ all lie in~$\mathfrak c$ and $\mathfrak c$~is
  $\sigma$-invariant, so we can find $Y \in \mathfrak k_s^{\mathbb C}$ such
  that $\langle X_n, Y, [X_n,Y] \rangle$ is an $\Sl(2,\mathbb C)$-triple
  commuting with~$\mathfrak b$.  As in the proof of
  Proposition~\ref{prop:nil-mon} we may find a real structure~$\sigma_s^g$
  on~$\mathfrak k_s^{\mathbb C}$ giving a compact form~$K_s^g$ and such
  that $Y=-\sigma_s^g X_n$ and $\sigma_s^g\mathfrak b= \mathfrak b$.  A
  transverse slice to~$\mathcal O_2$ at~$X$ is given by
  \begin{equation*}
    \begin{split}
      S_{X,Y}
      &= X + \mathfrak z(Y) \cap \mathfrak z(X_s) \\
      &= X_s + (X_n + \mathfrak z_{\mathfrak k_s^{\mathbb C}} (Y) ).
    \end{split}
  \end{equation*}
  Orbits of~$G$ meet~$S_{X,Y}$ in orbits of~$B$ and these lie in
  $K_s^g$-orbits.  The inequalities on cohomogeneities now follow as in
  Proposition~\ref{prop:nil-mon}. 
\end{proof}

\begin{proposition}
  \label{prop:mc2}
  If $G$~is a compact simple Lie group and $\mathcal O \subset \mathfrak
  g^{\mathbb C}$ is an orbit which is neither semi-simple nor nilpotent,
  then
  \begin{equation*}
    \cohom \mathcal O > 2.
  \end{equation*}
\end{proposition}

\begin{proof}
  Suppose $\cohom \mathcal O \leqslant 2$.  As $\mathcal O$ is not
  semi-simple, we have $\mathcal O \succneqq \mathcal O_s$ and
  Proposition~\ref{prop:mixed-mon} implies $\cohom \mathcal O > \cohom
  \mathcal O_s$.  Hence $\cohom \mathcal O_s$ is one.
  By~\cite{Dancer-Swann:hK-cohom1}, we have $G=\SU(n+1)$ and $\mathcal
  O_s$~is the orbit of $\diag(\lambda,\dots,\lambda,-n\lambda)$.  The mixed
  orbits with this semi-simple part are parameterised by nilpotent orbits
  of~$\Sl(n,\mathbb C)$, and have the same partial ordering.  This follows
  directly from the Jordan normal form of such elements.
  Proposition~\ref{prop:mixed-mon} implies that it is sufficient to
  calculate the cohomogeneity of the $\SL(n+1,\mathbb C)$-orbit, or
  equivalently $\GL(n+1,\mathbb C)$-orbit, of
  \begin{equation*}
    X=
    \begin{pmatrix}
      \lambda & 1 & & & \\
      & \lambda & &\smash{\hbox to 0pt{\hss\LARGE0\hss}} & \\
      & & \ddots & & \\
      &\smash{\hbox{\hss\LARGE0\hss}} & & \lambda & \\
      & & & & -n\lambda
    \end{pmatrix}
    .
  \end{equation*}
  The stabiliser of~$X$ in~$\Un(n+1)$ has Lie algebra $\mathfrak k=\un(1)_+
  + \un(n-2) + \un(1)_-$.  The fibre of the normal bundle to the real orbit
  of~$X$ at~$X$ is
  \begin{multline*}
    (\ad X)(\Sl(n+1,\mathbb C))/(\ad X)(\su(n+1))\\
    \cong
    \left\{
      \begin{pmatrix}
        0 & a & & &\\
        & 0 & & \smash{\hbox to 0pt{\hss\LARGE0\hss}} & \\
        & & \ddots& &\\
        &\smash{\hbox to 0pt{\hss\LARGE0\hss}} & & &\\
        b_1 & b_2 & \dots & b_n & 0
      \end{pmatrix}
      : a\in\mathbb R, b_i\in\mathbb C
    \right\}
    .
  \end{multline*}
  As a representation of~$\mathfrak k$, this fibre is 
  \begin{equation*}
    \mathbb R+ 2[L_+\overline{L_-}]+[\Lambda^{1,0}\overline{L_-}],
  \end{equation*}
  where $\Lambda^{1,0}$~is the standard representation of~$\Un(n-2)$
  on~$\mathbb C^{n-2}$.  This representation of~$\mathfrak k$ has
  cohomogeneity~$5$.
\end{proof}

We therefore do not need to consider mixed orbits any further.

\subsection{Orbits of Semi-Simple Groups}
\label{sec:ss-G}
If $G$~is a compact semi-simple Lie group, then the Lie algebra of~$G$
splits as
\begin{equation*}
  \mathfrak g = \mathfrak g_1 \oplus \dots \oplus \mathfrak g_r
\end{equation*}
with each $\mathfrak g_i$ simple.  Let $G_i$~be the simply-connected group
with Lie algebra~$\mathfrak g_i$.  If $\mathcal O$ is an adjoint orbit
for~$G^{\mathbb C}$, then we can write
\begin{equation*}
  \mathcal O = \mathcal O_1\times\dots\times\mathcal O_r,
\end{equation*}
where $\mathcal O_i$ is an adjoint orbit for~$G_i^{\mathbb C}$.  Now
\begin{equation*}
  \cohom_G\mathcal O = \cohom_{G_1}\mathcal O_1 +\dots+
  \cohom_{G_r}\mathcal O_r,
\end{equation*}
so we have the following result:

\begin{proposition}
  \label{prop:ssG-c2}
  If $G$~is a compact semi-simple group acting almost effectively on an adjoint
  orbit~$\mathcal O$ with cohomogeneity at most two, then either $G$~is
  simple and the orbit $\mathcal O$ is listed in \S\S\ref{sec:ss}
  and~\ref{sec:nil} or $G$~has two simple factors $G_1$ and~$G_2$, $\mathcal
  O=\mathcal O_1\times \mathcal O_2$ and $\mathcal O_i$ is of cohomogeneity
  one with respect to~$G_i$, for $i=1,2$.
  \qed
\end{proposition}

So in the semi-simple, non-simple case, the new orbits we must study are
$\mathcal O_1 \times \mathcal O_2$ where each $\mathcal O_i$ is either
the semi-simple orbit  $T^* \mathbb C \mathbb P(n)$ of $SL(n+1, \mathbb C)$
or a minimal nilpotent orbit.

\section{Quaternionic K\"ahler Structures}
\label{sec:qK-class}
The results of \S\S\ref{sec:moment} and~\ref{sec:orbits} give us the
following classification result.  At this stage we are not assuming
completeness of the metric.

\begin{theorem}
  \label{thm:open-simple-qK}
  Let $(M,g)$ be a quaternionic K\"ahler manifold of positive scalar
  curvature, of cohomogeneity one with respect to a compact simple
  group~$G$ acting effectively.
  
  Then up to covers, one of the following three alternatives holds\textup:
  \begin{enumerate}
  \item [(i)] $(M,g)$~is hyperHermitian, and $g$~is one of the metrics of
    \cite{Swann:MathAnn} as described in
    Proposition~\ref{prop:associated}\textup;
  \item [(ii)] an open set of the twistor space of~$(M,g)$ is an open set
    in a projectivised nilpotent orbit~$\mathbb P(\mathcal O)$ as a complex
    contact manifold with real structure and $\mathcal O$ is of
    cohomogeneity at most~$5$\textup;
    
    if the twistor space covers the whole of~$\mathbb P(\mathcal O)$, then
    $(M,g)$~is locally isometric to~$\mathfrak M(\mathcal O)$, and these
    cases are listed in Theorem~\ref{thm:nqK-c1}\textup;
  \item [(iii)] we are in the non-proper case, and the principal orbit
    of~$G$ in~$M$ is either
    \begin{gather*}
      \text{\textup{(a)}}\ \frac{\SU(n+1)}{\Un(n-1)},\quad
      \text{\textup{(b)}}\ \frac{\SU(n+2)}{\Un(n-2)}, \\
      \text{\textup{(c)}}\ \frac{\SO(10)}{\SU(2) \times \SU(2)},
      \quad\text{or}\quad
      \text{\textup{(d)}}\ \frac{\LieE_6}{\SU(4)}.
    \end{gather*}
  \end{enumerate}
\end{theorem}

\begin{proof}
  Part~(i), is the case of no open orbits on~$Z$ given in
  Proposition~\ref{prop:associated}.  Part~(ii) is the case of $Z$ having
  an open $G^{\mathbb C}$-orbit which is a proper complex contact
  manifold.  The statement about the real structure is due to Nitta \&\
  Takeuchi \cite{Nitta-Takeuchi:contact}.  When $Z$ covers all of~$\mathbb
  P(\mathcal O)$, we have that $Z$~is homogeneous under~$G^{\mathbb C}$ and
  the fact that $M$ is isometric to~$\mathfrak M(\mathcal O)$ was proved
  in~\cite{Swann:HTwNil}. 
  
  It only remains to deal with case~(iii).  Recall from~\S\ref{sec:moment}
  that in the non-proper case, we have a map $\mathbb Pf\colon U_1 \to
  \mathbb P (\mathcal O)$ with fibre a non-compact complex Lie group of
  complex dimension one.  We derived the inequalities
  \begin{equation}
    3 \geqslant \cohom_G(U_1) \geqslant \cohom_G(\mathbb P(\mathcal O))+1=
    \cohom_G(\mathcal O)+1.
    \label{eq:ss-3}
  \end{equation}
  By Proposition~\ref{prop:mc2}, the orbit~$\mathcal O$ is semi-simple, and
  hence the principal orbit types on~$\mathcal O$ are given by
  \eqref{eq:ss-p-orbit-1} and~\eqref{eq:ss-p-orbit-2}.
  
  Consider the principal orbit types $\SP(n+1)/\SP(n-1)$ and
  $\SO(n+2)/\SO(n-2)$, with $n \neq 4$ in the second case.  In these cases,
  the stabiliser acts trivially on the fibre of~$\mathbb Pf$, which has
  complex dimension one.  So the second inequality in~\eqref{eq:ss-3} is
  strict, and we get a contradiction.  Thus we are left with the cases
  listed in the Theorem (up to covers, the orthogonal case with $n=4$ gives
  the same principal orbit type in $M$ as case (b) with $n=2$).
\end{proof}

\begin{remark}
  We do not know whether all the cases in the above Theorem can be
  realised.  For case~(iii), we will see below that there is a compact
  example with special orbit as in (iii)(a).  We do not know whether
  (iii)(b--d) occur in non-compact examples.  In case~(ii), it is
  conceivable that open sets of the projectivised nilpotent orbits admit
  families of twistor lines not corresponding to the structures~$\mathfrak
  M(\mathcal O)$.  This is known to happen in~$\so(4,\mathbb C)$
  (see~\cite{Swann:HTwNil}), but currently we have no examples with
  $\mathfrak g^{\mathbb C}$ simple.
\end{remark}

Let us now consider the case when $M$~is compact.  For the nilpotent case
we will need an interpretation of a result of Beauville which is of
independent interest.  Note this result does not assume that $M$ is of
cohomogeneity one.

\begin{theorem}
  \label{thm:nilp-compact}
  Suppose $M$~is a compact quaternionic K\"ahler manifold with an action of a
  compact connected Lie group~$G$ preserving this structure.  If
  $G^{\mathbb C}$~has an open orbit on~$Z$ containing a twistor line then
  $M$~is isometric to a Wolf space and either $G$~is the isometry group of
  $M$ or $G$~is given by one of the shared-orbit pairs of Brylinski \&\ 
  Kostant \cite{Brylinski-Kostant:nilpotent-announcement,%
  Brylinski-Kostant:nilpotent}.
\end{theorem}

\begin{proof}
  The hypotheses on the action mean that $M$ is quaternionic K\"ahler with
  positive scalar curvature, so the twistor space~$Z$ is projective with
  $c_1(Z) > 0$, and $G^{\mathbb C}$ acts on~$Z$ algebraically.  The
  assumption that there is an open $G^{\mathbb C}$-orbit containing a
  twistor line implies that the image of~$\mathbb Pf$ is the
  projectivisation of a nilpotent orbit, i.e., we are in the proper case
  of~\S\ref{sec:moment}.  In particular, $\mathbb Pf$~has injective
  differential on an open set of~$Z$.  By \cite[III.5 Corollary~2,
  p.~252]{Mumford:red} the algebraic map~$\mathbb Pf$ is \'etale and
  hence finite.
  
  Exactly the same remarks apply to the component of the identity~$G_1$ of
  the full isometry group of~$M$ and the associated map $\mathbb Pf_1\colon
  Z\dashrightarrow \mathbb P(\mathfrak g_1^{\mathbb C})$.  Beauville's
  results~\cite{Beauville:Fano} apply to this $G_1^{\mathbb C}$-action
  on~$Z$, giving that $Z$~is the projectivised minimal nilpotent
  orbit~$\mathbb P(\mathcal O_1)$ for~$G_1^{\mathbb C}$ and that $M$~is the
  Wolf space with isometry group~$G_1$.
  
  If $G=G_1$ we are finished.  Otherwise $\mathbb Pf$~induces a finite
  morphism $\phi\colon\mathcal O_1 \dashrightarrow \mathcal O$ which is
  $G^{\mathbb C}$-equivariant.  This is the condition that $\mathcal O$~be
  a ``shared orbit'' and Beauville shows that only the pairs  listed
  in~\cite{Brylinski-Kostant:nilpotent-announcement,%
  Brylinski-Kostant:nilpotent} arise.
\end{proof}

\begin{theorem}
  Let $(M,g)$ be a compact quaternionic K\"ahler manifold, of cohomogeneity
  one with respect to a compact simple group~$G$.
  
  Then $(M,g)$ is a quaternionic K\"ahler symmetric space and the pair
  $(M,G)$ is given in Table~\ref{tab:compact-simple-qK}.
\end{theorem}

\begin{table}[htbp]
  \begin{center}
    \leavevmode
    \begin{tabular}{|c|c||c|c|}
      \hline
      $M$&$G$&$M$&$G$\\
      \hline
      \hline
      \vrule height 14pt width 0pt
      $\HP(n)$&$\SP(n)$&$\Gro_4(\mathbb R^7)$&$\LieG_2$ \\
      \vrule height 12pt width 0pt
      $\HP(n)$&$\SU(n+1)$&$\LieG_2/\SO(4)$&$\SU(3)$ \\ 
      \vrule height 12pt width 0pt
      $\Gr_2(\mathbb C^n)$&$\SU(n-1)$&$\LieF_4/\SP(3)\SP(1)$&$\Spin(9)$ \\
      \vrule height 12pt width 0pt
      $\Gr_2(\mathbb C^{2n})$&$\SP(n)$&$\LieE_6/\SU(6)\SP(1)$&$\LieF_4$\\
      \vrule height 12pt width 0pt
      $\Gro_4(\mathbb R^n)$&$\SO(n-1)$&&\\
      \hline
    \end{tabular}
    \medskip
    \caption{Compact quaternionic K\"ahler manifolds~$M$ which are
    cohomogeneity one under an action of a simple group~$G$.}
    \label{tab:compact-simple-qK}
  \end{center}
\end{table}

\begin{proof}
  We consider the three cases of Theorem~\ref{thm:open-simple-qK}.
  Case~(i), when there is no open $G^{\mathbb C}$-orbit in $Z$, was dealt
  with in Theorem~\ref{thm:hc-compact} and gives the first entry in the
  table.
  
  If there is an open $G^{\mathbb C}$-orbit $U_1$, then we have $U_1 =
  Z\setminus D$ for some subvariety~$D$ of~$Z$.  We are either in the
  proper case~(ii) or the non-proper case~(iii).
  
  In~(ii), we can apply the proof of Theorem~\ref{thm:nilp-compact}.  This
  together with Theorem~\ref{thm:open-simple-qK} implies that in order to
  get quaternionic K\"ahler manifolds of cohomogeneity one, we need to
  consider shared orbit pairs $(\mathcal O,\mathcal O_1)\subset(\mathfrak
  g^{\mathbb C},\mathfrak g_1^{\mathbb C})$ with $\mathcal O_1$ a minimal
  orbit and $\mathcal O$ next-to-minimal.  Looking at the list of shared
  orbits in~\cite{Brylinski-Kostant:nilpotent-announcement} we see that all
  but two cases give next-to-minimal orbits; the exceptions being the
  four-to-one covering of an orbit in $\so(8,\mathbb C)$ by the minimal
  orbit in $\lieF_4^{\mathbb C}$ and the six-to-one covering of the
  subregular orbit in~$\lieG_2^{\mathbb C}$ by the minimal orbit
  of~$\so(8,\mathbb C)$.  This gives the last six cases in the table.
  
  It remains to discuss case~(iii) and account for the second and third
  entries of the table, which we do using an argument of Poon \&\
  Salamon~\cite{Poon-Salamon:8}.  As mentioned above, the twistor space~$Z$
  is projective algebraic, with a linear action of~$G^{\mathbb C}$.
  Applying the Borel fixed point theorem to the action of a maximal
  torus~$T$ of~$G$, we see that $T$ has a fixed point on~$Z$, and hence
  on~$M$. Therefore the rank of the stabiliser of some point of~$M$ equals
  the rank of~$G$.
  
  As in the proof of Theorem~\ref{thm:hc-compact}, we see that $M/G$~is a
  compact interval, that there are two special orbits $G/H_1$ and~$G/H_2$,
  and, writing $G/K$ for the principal orbit, that $H_i/K$~is a sphere for
  $i=1,2$.  Examining the ranks of isotropy groups we see that cases
  (iii)(b), (c) and~(d) cannot arise.
  
  The remaining case is when $G=\SU(n+1)$ and the principal stabiliser has
  Lie algebra~$\un(n-1)$.  As the classification of compact quaternionic
  K\"ahler manifolds is known in dimensions four \cite{Hitchin:Kaehlerian}
  and eight \cite{Poon-Salamon:8}, we take $n>2$.
  
  Let us first assume that the principal stabiliser~$K$ is connected.  The
  long exact homotopy sequence for the fibration $K \rightarrow G
  \rightarrow G/K$, and the Hurewicz theorem, show that the second Betti
  number of~$G/K$ is~$1$.
  
  If $G/H$ is a special orbit then $H/K$ is a sphere~$S^r$ with~$r>0$,
  since $M$ is simply-connected.  Hence $H$ is connected, and the homotopy
  sequence for $H \rightarrow G \rightarrow G/H$ shows that $b_2(G/H)$ is
  the rank of $\pi_2(G/H) \cong \pi_1(H)$.  The latter may be computed by
  considering the fibration $K \rightarrow H \rightarrow S^r$ and we get
  that
  \begin{equation*}
   b_2(G/H) = 
   \begin{cases}
     2,&\text{if $r=1$},\\
     0,&\text{if $r=2$},\\
     1,&\text{if $r>2$}.
   \end{cases}
  \end{equation*}
  
  Following \cite{Alekseevsky-Podesta:co1}, we consider the Mayer-Vietoris
  sequence for an open cover of~$M$ consisting of the complements of each
  special orbit.  We obtain
  \begin{equation*}
    H^2(M) \rightarrow H^2\left(\frac G{H_1}\right) \oplus H^2 \left(\frac
      G{H_2}\right) \rightarrow H^2\left(\frac GK\right) \rightarrow H^3(M)
    = 0
  \end{equation*}
  in real cohomology.  If $M$ is not the complex Grassmannian then,
  by~\cite{LeBrun-Salamon:rigidity}, $b_2(M)=0$ and so one of the special
  orbits has $r=2$ and the other has $r>2$.
  
  Now, the complex coadjoint orbit $T^*\CP(n)$ has a special
  $\SU(n+1)$-orbit of real dimension~$2n$, so $Z$~contains a special
  $\SU(n+1)$-orbit of real dimension less than or equal to $2n+1$. It
  follows that one of the special orbits in~$M$ has real codimension at
  least $2n-1$.
  
  Therefore if $n>2$ and $M$ is not the complex Grassmannian then the
  special orbits in~$M$ must be $\SU(n+1)/\Un(n)$ and $\SU(n+1)/\SU(2)
  \times \SU(n-1)$ and $M$ is~$\HP(n)$ \cite{Alekseevsky-Podesta:co1}.
  
  If $K$ is disconnected we can argue as follows. The homotopy sequence for
  $K \rightarrow H \rightarrow S^r$ shows that either $r=1$ or $H$ is
  disconnected.  But if one special orbit $G/H_1$ has $r > 1$ then the
  codimension of that orbit is strictly greater than~$2$ and so the other
  special orbit~$G/H_2$, like $M$, must be simply-connected.  In
  particular, $G/H_2$ would have connected stabiliser and $r=1$.  The
  arguments of the preceding paragraphs now show that $b_2(G/H_2) = 2$ and
  $b_2(G/K) \leqslant 1$, so the Mayer-Vietoris argument again forces
  $M$~to be the complex Grassmannian.
  
  The remaining possibility, that $r=1$ for both special orbits, is ruled
  out for $n> 1$ because we know one of the special orbits in $M$ has
  codimension at least $2n-1$.
\end{proof}

\begin{remark}
  If $n=1$ then disconnected stabilisers can arise. We can view $S^4$ as an
  $\SO(3)$-manifold with both special orbits being~$\RP(2)$.  Also $\CP(2)$
  admits an $\SO(3)$-action with one orbit an~$S^2$ and the other
  an~$\RP(2)$.
\end{remark}

Let us now discuss the case when $G$~is semi-simple and $M$~is compact.  

\begin{theorem}
  Suppose $M$~is a compact quaternionic K\"ahler manifold of cohomogeneity
  one with respect to a compact group~$G$.  If $G$ is semi-simple but not
  simple, then $M=\HP(N)$ and $\mathfrak g = \LSP(N-m)\oplus\LSP(m+1)$ for
  some $0\leqslant m<N$.
\end{theorem}

\begin{proof}
  The case when there is no open $G^{\mathbb C}$-orbit on~$Z$ is covered by
  the results of~\S\ref{sec:no-open} and yields no examples with $G$~not
  simple.
  
  In the proper case, we may apply Theorem~\ref{thm:nilp-compact}.  The
  only shared orbit pair for a non-simple group occurs when
  $G=\SP(n_1)\times\dots\times\SP(n_r)$ and $\mathcal O$~is the product of
  the minimal orbits in $\LSP(n_1,\mathbb C),\dots,\LSP(n_r,\mathbb C)$.
  This is covered by the minimal orbit in~$\LSP(n_1+\dots +n_r,\mathbb C)$
  and $M$~is the quaternionic projective space~$\HP(n_1+\dots+n_r-1)$.
  This is of cohomogeneity one only when $G$~has exactly two factors, i.e.,
  $r=2$.
  
  In the non-proper case, the new orbits we must consider are,
  from~\S\ref{sec:moment}, of the form $\mathcal O_1 \times \mathcal O_2$
  where each of $\mathcal O_i$ is either the $\SL(n+1, \mathbb C)$-orbit
  $T^*\CP(n)$ or a minimal nilpotent orbit.  The principal orbit of
  $\SU(n+1)$ in the former orbit is $\SU(n+1)/\Un(n-1)$, while for the
  minimal nilpotent orbits for $H^{\mathbb C}$ the principal $H$-orbit is
  $H/K$, where $H/K\SP(1)$ is the associated Wolf space.
  
  In all these cases it follows that the rank of~$G$ is at least two
  greater than the rank of the principal stabiliser for the $G$-action on
  the orbit.  Also, the orbits are cohomogeneity two with respect to~$G$.
  Hence the rank of~$G$ is at least three greater than that of the
  principal stabiliser for the $G$-action on $Z$ or~$M$.  So the Borel
  fixed point argument shows that we do not obtain any compact examples
  with $G$~not simple.
\end{proof}

\begin{remark}
  Attempting to prove the above theorem for general reductive~$G$ requires
  stronger results than we were able to obtain in~\S\ref{sec:no-open}.
\end{remark}

\section{HyperK\"ahler Structures}
\label{sec:hK-class}
The results obtained in~\S\ref{sec:orbits} can be used to make results of
Bielawski~\cite{Bielawski:homogeneous} a little more precise.  In
\cite{Dancer-Swann:hK-cohom1} we were able to classify hyperK\"ahler metrics
of cohomogeneity one.  We now wish to extend this to metrics of
cohomogeneity two.

Suppose $M$~is a hyperK\"ahler manifold with Riemannian metric~$g$ and
complex structures $I$, $J$ and~$K$ satisfying $IJ=K=-KI$.  Then
$aI+bJ+cK$~is also a complex structure on~$M$ whenever $(a,b,c)\in S^2\subset
\mathbb R^3$.  Let $\mathbf S$ denote the set of these complex structures.
Suppose that a compact semi-simple Lie group~$G$ acts on~$M$ preserving
both the metric and the individual complex structures and that the action
is of cohomogeneity two.  Fix a point~$x$ of~$M$ lying on some principal
orbit.  Then at~$x$ the tangent vectors generated by the $G$-action span a
real subspace~$\mathcal D$ of~$T_xM$ of codimension~$2$.

We claim that there are at most two points $\mathcal J$, $-\mathcal J$
in~$\mathbf S$ which preserve~$\mathcal D$.  Suppose $A, B\in \mathbf S$
are linearly independent and preserve~$\mathcal D$.  Then $\lambda A+\mu B$
also preserves~$\mathcal D$ for any $\lambda, \mu\in\mathbb R$, so we may
assume $A$ and~$B$ are orthogonal.  But now $AB=-BA$ and this
preserves~$\mathcal D$, so $\mathcal D$~is an $\mathbb H$-module,
contradicting the fact that $\mathcal D$~has real codimension~$2$ in the
$\mathbb H$-module~$T_xM$.

Thus, without loss of generality, we may assume that $I$ does not
preserve~$\mathcal D$.  Infinitesimally, we may complexify the action of~$G$
with respect to~$I$.  Then the tangent vectors to this action span~$T_xM$
and they also span the tangent space of~$M$ in a $G$-invariant
neighbourhood of~$x$.  The assumption that $G$~be semi-simple implies that
we have a hyperK\"ahler moment map 
\begin{equation*}
  \mu = \mu_I i + \mu_J j + \mu_K k \colon M \to \mathfrak g\otimes
  \im\mathbb H,
\end{equation*}
where we have identified $\mathfrak g^*$ with $\mathfrak g$ via the Killing
form.  The map $\mu_c:=\mu_J+i\mu_K$ is a moment map for the
complex-symplectic action of~$G^{\mathbb C}$.  In a neighbourhood of~$x$,
the map $\mu_c$~is a local diffeomorphism to an open subset of an adjoint
orbit $\mathcal O\subset \mathfrak g^{\mathbb C}$.  Moreover, $\mu_c$~is
$G$-equivariant, so the orbit~$\mathcal O$ is of cohomogeneity two and is
given by Proposition~\ref{prop:ssG-c2}.  Note that each of these orbits
actually arises, since Kronheimer, Biquard and Kovalev have shown that
every adjoint orbit carries a $G$-invariant hyperK\"ahler metric compatible
with the Kirillov-Kostant-Souriau complex-symplectic form
\cite{Kronheimer:nilpotent,Kronheimer:semi-simple,Biquard:orbits,Kovalev:Nahm}.

If we assume that $M$~is complete and that for some choice of~$I$ the
$G^{\mathbb C}$-action is locally transitive, then
Bielawski~\cite{Bielawski:homogeneous} shows that $\mathcal O$~is a
semi-simple orbit and that for one of the hyperK\"ahler metrics on~$\mathcal
O$ constructed by Kronheimer~\cite{Kronheimer:semi-simple}, $\mu_c\colon
M\to\mathcal O$~is an isometry.  

Noting that the above discussion also applies to $G$ with central factors
providing we assume the existence of a hyperK\"ahler moment map, we have the
following result.

\begin{theorem}
  Suppose $M$~is a complete hyperK\"ahler manifold which is of
  cohomogeneity two with respect to an isometric action of a compact
  group~$G$ preserving each complex structure.  If $G$~is semi-simple or,
  more generally, if there is a hyperK\"ahler moment map for the action
  of~$G$ on~$M$, and if for some~$I$ the $G^{\mathbb C}$-action is
  locally transitive, then $M$~is isometric to one of\textup:
  \begin{enumerate}
  \item [(1)] $T^*\CP(2n-1)$, with $G=\SP(n)$\textup;
  \item [(2)] $T^*\Gr_2(\mathbb C^n)$, with $G=\SU(n)$\textup;
  \item [(3)] $T^*\Gro_2(\mathbb R^n)$, with $G=\SO(n)$\textup;
  \item [(4)] $T^*\bigl(SO(10)/U(5) \bigr)$\textup;
  \item [(5)] $T^*\bigl(\LieE_6/\Spin(10)\SO(2)\bigr)$\textup;
  \item [(6)] $T^*\CP(n)\times T^*\CP(m)$, with $G=\SU(n+1)\times\SU(m+1)$.
  \end{enumerate}
\end{theorem}

\section{Three-Sasakian Structures}
\label{sec:3S-class}
Let $S$ be a compact $3$-Sasakian manifold with Riemannian metric~$g$.  We
refer the reader to~\cite{Boyer-GM:three-Sasakian} for a precise
definition, but for our purposes the following characteristic property will
be sufficient: the manifold $N_S:= S\times \mathbb R_{>0}$ with the metric
$dr^2+r^2g$ is hyperK\"ahler with an action of~$\SP(1)$, trivial on the
$\mathbb R_{>0}$-factor and satisfying (i)--(iv) of
Proposition~\ref{prop:sp1}.  This gives examples of $3$-Sasakian manifolds
by taking $N_S=\UM(M)$ with $M$ a quaternionic K\"ahler manifold of positive
scalar curvature.  In general, $S/\SP(1)$~is only a quaternionic K\"ahler
orbifold even though $S$~is smooth.

The group of $3$-Sasakian symmetries of~$S$ is exactly the group of
triholomorphic isometries of~$N_S$ and is necessarily compact.  Suppose $G$
is a compact group of symmetries of~$S$.  Then there is a unique
hyperK\"ahler moment map $\mu\colon N_S \to \mathfrak g^* \otimes \im\mathbb
H$ for this action with the property that
\begin{equation}
  \label{eq:scaling}
  \mu(s,r)=r^2\mu(s,1),
\end{equation}
for $s\in S$, $r\in\mathbb R_{>0}$ \cite{Boyer-GM:three-Sasakian}.  Recall
that $\mu$ is $G$-equivariant.

Suppose $S$ is of cohomogeneity one under the action of~$G$.  Then the
hyperK\"ahler manifold~$N_S$ is of cohomogeneity two and we may apply some of
the results of the previous section.

Let $P_S$ and $P_N$ denote the union of principal orbits of $S$ and~$N_S$
respectively.  Note that as $S$ is compact and Einstein with positive
scalar curvature, the fundamental group $\pi_1(S)$ is finite and so $P_S$
can not be all of~$S$.

\begin{lemma}
  \label{lem:S-I}
  The hyperK\"ahler manifold $N_S$ has a compatible complex structure~$I$
  such that $T_xN_S = T(G\cdot x) + IT(G\cdot x)$ for all $x\in P_N$.
\end{lemma}

\begin{proof}
  In \S\ref{sec:hK-class} we showed that at each $x\in P_N$ the set
  \begin{equation*}
    \mathcal A_x=\left\{\,\mathcal J\in \mathbf S: T_xN_S\ne T(G\cdot
      x)+\mathcal JT(G\cdot x)\,\right\}
  \end{equation*}
  is either empty or a point of~$\RP(2)=\mathbf S/\{\pm1\}$.  However,
  $\mathcal A_{(s,r)}=\mathcal A_{(s,1)}$ for $r\in\mathbb R_{>0}$ and
  $\mathcal A_{(g\cdot s,r)}=\mathcal A_{(s,r)}$ for $g\in G$, so we have a
  map from a subset of $P_S/G=P_N/(G\times \mathbb R_{>0})\to \RP(2)$.
  This map is a smooth map of a one-dimensional manifold, so its image is
  one-dimensional and the set $\mathbf S\setminus \cup_{x\in P_N} \mathcal
  A_x$ of $I$ satisfying our requirements is non-empty.
\end{proof}

Fix such a choice of~$I$ and let $J,K\in\mathbf S$ be such that $IJ=K=-JI$.
Let $\mu_I$ be the component of the hyperK\"ahler moment map~$\mu$
corresponding to~$I$, etc.  The group $\mathbb C^*$ acts on~$N_S$ as the
complexification with respect to~$I$ of the $\mathbb R_{>0}$-action.  The
quotient $Z=N_S/\mathbb C^*$ is the \emph{twistor space} of~$S$ and is a
compact, normal, $\mathbb Q$-factorial Fano variety~\cite{Boyer-G:twistor}
and, in particular, is projective.  The hyperK\"ahler manifold $N_S$ is the
complement of the zero section of an algebraic line bundle on~$Z$.  We now
have that $G^{\mathbb C}$ acts on both $Z$ and~$N_S$.

As $P_N$ is connected, Lemma~\ref{lem:S-I} implies that there is an open
$G^{\mathbb C}$-orbit~$U$ containing~$P_N$.  The scaling
property~\eqref{eq:scaling} implies that $\mu_c=\mu_J+i\mu_K$ maps~$U$ to a
nilpotent orbit~$\mathcal O$.  Note that $\mu_c$ is $\mathbb
C^*$-equivariant.  The usual rank argument for the moment maps shows that
the fibres of~$\mu_c$ on~$U$ are of dimension zero.  In particular, $U$~and
$\mathcal O$ have the same dimension.  Any other $G^{\mathbb C}$-orbit
on~$N_S$ lies in the complement of~$P_N$, so has smaller dimension
than~$\mathcal O$, and hence can not lie in~$\mu_c^{-1}(\mathcal O)$.  Thus
$U=\mu_c^{-1}(\mathcal O)$.

\begin{lemma}
  \label{lem:finite-g}
  The map $\mu_c$~is an algebraic morphism, finite over~$\mathcal O$.  The
  symmetry group~$G$ is semi-simple with at most two simple factors.  If
  $G$~is simple then $\mathcal O$ is next-to-minimal and $\mathfrak
  g\ne\su(3)$.  If $G$~is not simple then $\mathcal O$ is the product of a
  minimal orbit in each simple factor of~$\mathfrak g^{\mathbb C}$.
\end{lemma}

\begin{proof}
  The definition of~$\mu_c$ and the algebraic nature of~$N_S$ imply that
  $\mu_c$~is an algebraic map and so, by the above discussion, is
  finite-to-one on~$\mu_c^{-1}(\mathcal O)$.  Thus only the semi-simple
  part of the component of the identity of~$G$ acts effectively.  As
  $\mathcal O$ has cohomogeneity two with respect to~$G$, the list of
  possible nilpotent orbits~$\mathcal O$ follows from the results
  of~\S\ref{sec:orbits}.
\end{proof}  

Let $G_1$ be the component of the identity of the full 3-Sasakian
 symmetry group of~$S$.  Then $G_1$ satisfies the same hypotheses
 as~$G$, but might act homogeneously on~$S$.  Let $\nu\colon N_S\to
 \mathfrak g_1^*\otimes\im H$ be the moment map for the action
 of~$G_1$.

\begin{lemma}
  \label{lem:finite-1}
  The map $\nu_c$ is a finite map without zeroes.  The image of~$\nu_c$ is
  $\overline{\mathcal O_1}\setminus\{0\}$, where $\mathcal O_1$ is a
  nilpotent orbit in~$\mathfrak g_1^{\mathbb C}$ which is either minimal or
  listed in Lemma~\ref{lem:finite-g}.
\end{lemma}

\begin{proof}
  If $G_1$~acts transitively on~$S$ then there is nothing to prove as
  $G_1^{\mathbb C}$ has only one orbit on~$N_S$.  So we suppose that
  $G_1$ acts with cohomogeneity one on~$S$.
  
  Let $Q=N_S \setminus \nu_c^{-1}(\mathcal O_1)$.  Then $Q$~is complex, as
  it is the pre-image of the subvariety $\overline{\mathcal O_1} \setminus
  \mathcal O_1$ of~$\overline{\mathcal O_1}$.  As $\nu_c^{-1}(\mathcal
  O_1)$ contains the set of principal $G_1$-orbits and $\nu_c$~has the
  scaling property~\eqref{eq:scaling}, the set~$Q$ is a union of components
  of $(G_1/H_1 \coprod G_1/H_2) \times \mathbb R_{>0}$, where $G_1/H_i$ are
  the two special orbits in~$S$.  As $G_1$~commutes with the action
  of~$\SP(1)$, we conclude that the decomposition
  \begin{equation}
    \label{eq:N-Q}
    N_S= \nu_c^{-1}(\mathcal O_1) \cup Q
  \end{equation}
  is invariant under the action of~$\SP(1)\times\mathbb R_{>0}$.  In
  particular, the decomposition~\eqref{eq:N-Q} is determined by the action
  of the real group~$G_1$ and not by the particular choice of complex
  structure~$I$.  So $Q$~is complex for $I$, $J$ and~$K$ and the real
  dimension of~$Q$ is divisible by four.  Note that $G_1$~acts on~$Q$ with
  cohomogeneity one, so each component of~$Q$ is a single $G_1^{\mathbb
  C}$-orbit.

  The scaling property~\eqref{eq:scaling} implies that $\nu_c$ induces a
  map 
  \begin{equation*}
    \mathbb P\nu\colon Z\dashrightarrow \mathbb P(\mathfrak g^{\mathbb
    C}_1).
  \end{equation*}
  The set where $\mathbb P\nu$ fails to be defined lies in~$\mathbb PQ$ and
  so has complex codimension~$2$ in~$Z$.  The hypotheses of
  \cite[Lemma~3.3]{Beauville:Fano} now apply to~$\mathbb P\nu$.  We
  conclude that $\mathbb P\nu$~is everywhere defined and that both $\mathbb
  P\nu$ and~$\nu_c$ are finite.  This implies, using the scaling property,
  that $\nu_c$~has no zeroes.  The image of~$\mathbb P\nu$ is the closure
  of~$\mathbb P\mathcal O_1$ and so the image of~$\nu_c$ is
  $\overline{\mathcal O_1} \setminus \{0\}$.  The orbit $\mathcal O_1$ has
  cohomogeneity at most two with respect to~$G_1$ and so is either minimal
  or next-to-minimal.
\end{proof}

\begin{theorem}
  Suppose $S$ is a compact $3$-Sasakian manifold which is of cohomogeneity
  one under a compact connected Lie group~$G$ preserving the $3$-Sasakian
  structure.  Then $S$~is homogeneous under its full $3$-Sasakian symmetry
  group and the possible pairs $(S,G)$ are given in Table~\ref{tab:co1-3S}.
\end{theorem}

\begin{table}[htbp]
  \begin{center}
    \begin{tabular}{|c|c|}
      \hline \vrule height 12pt width 0pt
      $S$&$G$\\
      \hline \hline \vrule height 14pt width 0pt
      $S^{4n+3}$ & $\SP(r)\times\SP(n+1-r)$  \\
      \vrule height 12pt width 0pt
      $\RP(4n+3)$ & $\SP(r)\times\SP(n+1-r)$  \\
      \vrule height 12pt width 0pt
      $\SO(n+1)/\SO(n-3)\SP(1)$ & $\SO(n)$  \\
      \vrule height 12pt width 0pt
      $\SU(2n)/\Special(\Un(2n-2)\Un(1))$ & $\SP(n)$  \\
      \vrule height 12pt width 0pt
      $\SO(7)/\SO(4)\SP(1)$ & $G_2$  \\
      \vrule height 12pt width 0pt
      $\LieF_4/\SP(3)$ & $\Spin(9)$  \\
      \vrule height 12pt depth 4pt width 0pt
      $\LieE_6/\SU(6)$ & $\LieF_4$  \\
      \hline
    \end{tabular}
    \medskip
    \caption{Compact three-Sasakian manifolds~$S$ of cohomogeneity one
    under a connected compact Lie group~$G$ of $3$-Sasakian symmetries.}
    \label{tab:co1-3S}
  \end{center}
\end{table}

\begin{proof}
  Let $\mathcal O_1$ be as above.  Following Beauville, we will show that
  $\mathcal O_1$ is a minimal nilpotent orbit and that $(\mathcal
  O,\mathcal O_1)$ is a shared-orbit pair.  
  
  The space $N_S$ is a smooth finite branched $G_1^{\mathbb C}$-equivariant
  cover of ${\overline{\mathcal O_1}}^\times := \overline{\mathcal O_1}
  \setminus \{0\}$, such that $S=N_S/\mathbb R_{>0}$ is smooth, compact and
  cohomogeneity at most one under~$G_1$.  We determine all such covers for
  the orbits~$\mathcal O_1$ of Lemma~\ref{lem:finite-1}.  A priori this is
  more general than $\overline{\mathbb P\mathcal O_1}$ having a smooth
  equivariant cover; however, much of Beauville's paper applies at the
  level of orbits.
  
  Firstly, if $\nu_c\colon N_S \to {\overline{\mathcal O_1}}^\times$ is
  birational then the proof of Corollary~5.3 in~\cite{Beauville:Fano} shows
  that $\mathcal O_1$ is either the next-to-minimal orbit
  in~$\lieG_2^{\mathbb C}$ or it is a minimal nilpotent orbit.
  
  Now suppose that $\nu_c\colon N_S\to {\overline{\mathcal O_1}}^\times$
  has non-trivial degree~$\delta$.  Note that $\delta$
  divides~$\abs{\pi_1(\mathcal O_1)}$.
  
  Consider the case when $G_1$ is simple.  Lemma~\ref{lem:finite-1} implies
  that $\mathcal O_1$ is either the minimal orbit in~$\LSP(n,\mathbb C)$ or
  a next-to-minimal orbit in a simple Lie algebra other than~$\Sl(3,\mathbb
  C)$ (since $\LSP(n,\mathbb C)$ is the only simple Lie algebra whose
  minimal orbit is not simply connected).  Using the tables
  in~\cite{Collingwood-McGovern:nilpotent} to find $\abs{\pi_1(\mathcal
  O_1)}$, we see that $\mathcal O_1$~is either (i)~$\mathcal
  O_{(31^{n-3})}$ in~$\so(n,\mathbb C)$, or (ii)~$\mathcal O_{(2^41)}$
  in~$\so(9,\mathbb C)$, or the unique next-to-minimal orbit in
  (iii)~$\LSP(n,\mathbb C)$ or (iv)~$\lieF_4^{\mathbb C}$.  In each of
  these cases $\pi_1(\mathcal O_1)=\mathbb Z/2$ and thus $\delta=2$.  The
  manifold $N_S$ is now uniquely determined by~$\mathcal O_1$ as the
  complement of the vertex in the normalisation of~$\overline{\mathcal
  O_1}$ in the function field of the universal cover of~$\mathcal O_1$
  (cf.~\cite{Brylinski-Kostant:nilpotent}).  In the four cases to hand,
  this space is the minimal nilpotent orbit~$\widehat{\mathcal O} \subset
  \widehat{\mathfrak g}^{\mathbb C}$ where $\widehat{\mathfrak g}^{\mathbb
  C}$ is either (i)~$\so(n+1,\mathbb C)$, (ii)~$\lieF_4^{\mathbb C}$,
  (iii)~$\su(2n,\mathbb C)$ or (iv)~$\lieE_6^{\mathbb C}$.

  If $G_1$~is not simple, then $\mathcal O_1=\mathcal O_a\times\mathcal
  O_b$ with $\mathcal O_j$~the minimal orbit in the simple Lie
  algebra~$\mathfrak g_j^{\mathbb C}$, for $j=a,b$.  Dividing by the action
  of~$\mathbb R_{>0}$, we get a finite surjective $G_1$-equivariant map
  \begin{equation*}
    \mathbb R\nu\colon S\to S_1:= \left(\overline{(\mathcal
      O_a\times\mathcal O_b)} \setminus\{(0,0)\}\right) / \mathbb R_{>0}
  \end{equation*}
  induced by~$\nu_c$.  The special orbits of~$S_1$ have strictly smaller
  dimension than the principal orbits.  As $\mathbb R\nu$ is a finite map,
  we deduce that $\mathbb R\nu$ maps principal orbits to principal orbits
  and that $\mathbb R\nu$~is a bijection of the $G_1$-orbit spaces.  Now
  $S$~is smooth with principal orbits~$G_1/H$ and special orbits~$G_1/H_i$,
  so $H_i/H$ are spheres, $i=1,2$.  In $S_1$, we have principal orbits
  $G_a/K_a\times G_b/K_b$ and special orbits $G_a/K_a$ and $G_b/K_b$.  Thus
  one of the spheres $H_i/H$ finitely covers $(G_a\times K_b)/(K_a\times
  K_b) = G_a/K_a = \mathcal O_a/\mathbb R_{>0}$ in~$S_1$.  The other sphere
  covers $\mathcal O_b/\mathbb R_{>0}$.  So $\mathcal O_j/\mathbb R_{>0}$,
  $j=a,b$, is a finite quotient of a sphere, which is only the case for
  $\mathfrak g_j^{\mathbb C}=\LSP(n_j,\mathbb C)$.  As in
  \cite[Proposition~6.8]{Beauville:Fano} we now get that $N_S$~is the
  double cover of the minimal orbit of~$\LSP(n_a+n_b,\mathbb C)$.
  
  In summary, the above three cases show that $N_S$ is the cover of a
  minimal nilpotent orbit~$\widehat{\mathcal O}\subset \widehat{\mathfrak
  g}^{\mathbb C}$ with $\widehat{\mathfrak g}^{\mathbb C}$~simple: more
  precisely, either $N_S=\widehat{\mathcal O}$ or $N_S$ is the double cover
  of the minimal nilpotent orbit in~$\LSP(n,\mathbb C)$.  But now $Z$ is
  $\mathbb P\widehat{\mathcal O}$ as a complex contact manifold with real
  structure and so the $3$-Sasakian symmetry group of~$S$ is a compact real
  form~$\widehat G$ of~$\widehat{\mathfrak g}^{\mathbb C}$.  But $\widehat
  G$ acts transitively on $\widehat{\mathcal O}/\mathbb R_{>0}$.  So $S$~is
  homogeneous, $G_1=\widehat G$ and $\mathcal O_1=\widehat{\mathcal O}$.
  The pair $(\mathcal O,\mathcal O_1)$ is now a shared orbit pair, with
  $\mathcal O_1$ a minimal nilpotent orbit and $\mathcal O$ as in
  Lemma~\ref{lem:finite-g}.  Table~\ref{tab:co1-3S} now follows from the
  list in~\cite{Brylinski-Kostant:nilpotent-announcement}.
\end{proof}

\begin{remark}
  On a compact $3$-Sasakian manifold~$S$ it is also natural to consider the
  full group of isometries~$\Isom(S)$.  By results of
  Tanno~\cite{Tanno:Sasakian} (see~\cite{Boyer-GM:three-Sasakian}), if
  $S$~is not of constant curvature, the Lie algebra~$\mathfrak i$ of the
  isometry group is the direct sum of the Lie algebra~$\mathfrak s$ of the
  group~$\Sym(S)$ of $3$-Sasakian symmetries and $\LSP(1)$.  The
  $\SP(1)$-orbits are the fibres of the canonical fibration of~$S$ over a
  quaternionic K\"ahler orbifold~$M$.  In particular, the quotients of~$S$ by
  $\Isom(S)$ and of $M$ by~$\Sym(S)$ agree and we have
  \begin{equation*}
    \cohom_{\Isom(S)} S = \cohom_{\Sym(S)} M.
  \end{equation*}
  
  Working on unions of principal orbits we may apply some of our results
  when $S$ has a semi-simple group of isometries acting with cohomogeneity
  one.  In this case, there is a semi-simple subgroup~$G$ of~$\Sym(S)$ such
  that $G\times\SP(1)$ acts with cohomogeneity at most one on~$S$ and hence
  $\cohom_G M \leqslant 1$.  When $G^{\mathbb C}$ has no open orbit on the
  twistor space, the techniques of~\S\ref{sec:no-open} imply that there are
  no examples with non-constant curvature.  If $G^{\mathbb C}$ has an open
  orbit on~$Z$ then in the non-proper case we may apply the Borel fixed
  point arguments as in~\S\ref{sec:qK-class} and deduce that
  $\SU(n)\times\SU(2)$ acts with cohomogeneity one on~$S$.  Non-homogeneous
  examples occur in this case as $3$-Sasakian quotients of a sphere via a
  $\Un(1)$ action embedded in~$\SP(n+1)$ with weights $(1,\dots,1,k)$ as
  described in~\cite{Boyer-GM:three-Sasakian}.  We do not know whether
  these are the only such examples.  The proper case is the hardest to
  analyse.  We are unable to apply Beauville's techniques because there is
  no bound on the second Betti number of a $3$-Sasakian manifold and we
  have not been able to find other arguments to show that the complement of
  the open $G^{\mathbb C}$-orbit has codimension at least~$2$.
\end{remark}

\providecommand{\bysame}{\leavevmode\hbox to3em{\hrulefill}\thinspace}

\end{document}